\documentclass{article}

\usepackage{graphicx}
\usepackage{balance}  

\usepackage[utf8]{inputenc}
\usepackage{ dsfont }
\usepackage{ bbold }
\usepackage{ amsmath }
\usepackage{graphicx} 
\usepackage{caption}
\usepackage{subcaption} 

\usepackage{multirow}
\usepackage{mdframed} 

\usepackage[numbers]{natbib}

\allowdisplaybreaks 

\mdfdefinestyle{myenvs}{%
  hidealllines=false,
  nobreak=true, 
  leftmargin=0pt,
  rightmargin=0pt,
  innerleftmargin=0pt,
  innerrightmargin=0pt,
}

\newmdtheoremenv[style=myenvs]{theorem}{Theorem}

\newtheorem{bound}{Bound}

\begin{document}
\title{Q-error Bounds of Random Uniform Sampling for Cardinality Estimation}
\author{Beibin Li$^1$\footnote{The author performed the work while interning at Microsoft Research}, 
Yao Lu$^2$, Chi Wang$^2$, Srikanth Kandula$^2$}
\date{
$^1$University of Washington, 
$^2$Microsoft Research
}

\maketitle

\begin{abstract}
Random uniform sampling has been studied in various statistical tasks but few of them have covered the Q-error metric for cardinality estimation (CE). In this paper, we analyze the confidence intervals of random uniform sampling with and without replacement for single-table CE.  Results indicate that the upper Q-error bound depends on the sample size and true cardinality. Our bound gives a rule-of-thumb for how large a sample should be kept for single-table CE. 
\end{abstract}

\section{Introduction}
Cardinality estimation (CE) is the key to various tasks such as query optimization and approximate query processing. Production systems apply histograms~\cite{Ioannidis:2003:HH:1315451.1315455} and sketches such as Count-Min and HyperLogLog~\cite{Cormode:2005:IDS:1073713.1073718} for fast and accurate  estimates~\cite{kmvSynopsis,Cormode:2012:SMD:2344400.2344401,Bruno:2001:SMW:375663.375686,Cormode:2005:IDS:1073713.1073718,Flajolet:1985:PCA:5212.5215,numdistinct,hhstreaming}. 
Random uniform sampling has been studied in this context for a long period of time; despite various prior analyses for random sampling~\cite{Haas:1996:SCE:233502.233514,haas1995sampling}, few have covered the Q-error metric for CE.

In this paper, we analyze the Q-error bounds for sampling-based single-table CE with and without replacement.  Based on existing statistical tools such as the Chernoff Bound and the Bernstein-Serfling's Inequality, our analyses show the confidence intervals for the Q-error are less than a threshold given the sample size and true cardinality. The upper Q-error bound for random sampling with replacement is agnostic to the size of the original table. Our analyses can be used as a rule-of-thumb for how large a sample should be kept to reach a specific Q-error at a given confidence interval, as well as a simple accuracy baseline for other CE solutions~\cite{aidb}.

In the following paper, Section~\ref{sec:definition} defines the problem setup.  Section~\ref{sec:rs} and Section~\ref{sec:rs_np} analyze the upper error bounds of random  sampling with and without replacement for single-table CE. We discuss link to related work in Section~\ref{sec:related}.

\section{Problem Setup}
\label{sec:definition}
Let $t$ be a table with $n$ rows and $m$ columns, and $\hat t$ be a sample of $k$ rows drawn uniformly at random from $t$.
Suppose $pn$ rows from $t$ satisfy the predicate, where $p\in[0,1]$ represents the probability (aka selectivity) that a row in the table satisfies the given predicate.
Cardinality estimation (CE), i.e., estimating $np$ given a predicate, can be formulated as an application of Binomial distribution (Sum of Independent Bernoulli Trials), in which each random variable (corresponding to each row) takes the value of 1/0 for satisfying the predicate or not.
 
Let $X = pn = \sum\limits_{i=0}^n x_i$ be the cardinality of the predicate in the original table,  where $x_i = 1$ if row $i$  satisfies the predicate.  Let $\hat{X} =  \sum\limits_{i=0}^k \hat{x}_i$ be the cardinality of the predicate on the sampled table. We  use $\mu = \mathds{E}[\hat{X}] = kp$ as the expected number of rows that satisfy the predicate in the sampled table. The population variance is defined as $\sigma^2 = \frac{1}{n} \sum\limits_{i=0}^n (x_i - p)^2 = \mathds{E}[x_i^2] - \mathds{E}[x_i]^2 = (1^2 \times p + 0^2 \times (1-p)) - p^2= p (1-p) $. 
 
Following recent CE work~\cite{lm, naru, deepdb, mscn}, we use the Q-error metric in our analysis for evaluating the estimation accuracy:  
$${\tt Q\mbox{-}error}=\max(\frac{\texttt{true}}{\texttt{est}}, \frac{\texttt{est}}{\texttt{true}}),$$ 
where $\texttt{est}$ is the estimated cardinality and $\texttt{true}$ is the true cardinality. 
In practice, we replace $\texttt{est} = \max(\texttt{est}, 1)$ and $\texttt{true} = \max(\texttt{true}, 1)$ to avoid divide-by-zero. {\tt Q-error}=1 indicates a perfect prediction.

\section{Bounds of Sampling with Replacement}
\label{sec:rs}
 
Random uniform sampling with replacement under the independent and identically distributed (i.i.d.) assumption is widely applied in modern machine learning (e.g., bootstrap \citep{efron1986bootstrap}).
Applying Chernoff Bound to sampling with replacement gives us a concise bound with the Q-error metric. We also incorporate error bounds from Bernstein's Inequality to tighten the bound.

\begin{theorem}
\label{theory:with_replacement}
For cardinality estimation over single tables, the Q-error of random uniform sampling with replacement is bounded by
\begin{align*}
& \mathds{P}(\text{Q-error $\leq$ q}) \geq 1 - \Omega - \Psi, where \\
\Omega &=  \min\Bigg(\Big(\cfrac{e^{q - 1}}{q^{q}} \Big)^{pk},  \exp\Big(- \frac{k (pq-p)^2}{2 \sigma^2 + 2 (pq-p) / 3}\Big) \Bigg), \\
\Psi &=  \min\Bigg(\Big(e^{(\frac{1}{q} - 1)} q^{\frac{1}{q}} \Big)^{pk},  \exp\Big(- \frac{k (p-p/q)^2}{2 \sigma^2 + 2 (p-p/q) / 3}\Big) \Bigg).
\end{align*}
\end{theorem}
 
In the case when the probability $ \mathds{P}(\text{Q-error $\leq$ q}) $ is negative (i.e. $\Omega + \Psi > 1$), we replace it with zero.
$\Omega$ is the probability for over-estimation and $\Psi $ is the probability for under-estimation.
We further show the result with Hoeffding's inequality in Appendix, since adding it can only tighten the bounds slightly in some corner cases.
 
From Theorem \ref{theory:with_replacement}, we can see that the final error $q$ is \emph{agnostic} to the number of rows $n$ of the original table. The bound is only relevant to the number of samples $k$ and the true cardinality $p$. We give the detailed analyses below.

\vspace{0.1in}\noindent\textbf{CE bounded by Chernoff}. 
The Chernoff bound \cite[Reference][Corollary 4.2]{mulzer2018five} for the sum of Bernoulli trials can be stated as the following inequalities, where $\delta \in [0,1]$:
 
\begin{align}
\mathds{P}(\hat{X} \geq (1 + \delta) \mu) \leq \Big(\cfrac{e^\delta}{(1 + \delta)^{(1 + \delta)}} \Big)^\mu, \label{eq:cb:geq}\\
\mathds{P}(\hat{X} \leq (1 - \delta) \mu) \leq \Big(\cfrac{e^{-\delta}}{(1 - \delta)^{(1 - \delta)}} \Big)^\mu,  \label{eq:cb:leq}
\end{align} 
 
\noindent Let $q\ge 1$ be the q-error, and we have the over-estimation probability
 
\begin{align*}
\mathds{P}(Pred \geq q X) &=  \mathds{P}(\hat{X} \geq q \mu) \\
&\text{(Let $(1 + \delta) = q$ and apply (\ref{eq:cb:geq}))} \\
&\leq \Big(\cfrac{e^\delta}{(1 + \delta)^{(1 + \delta)}} \Big)^\mu \\
&= \Big(\cfrac{e^{q - 1}}{q^{q}} \Big)^{pk} = \Omega.
\end{align*}

\noindent Similarly, we have the under-estimation probability
\begin{align*}
\mathds{P}(Pred \leq \frac{1}{q} X) &=  \mathds{P}(\hat{X} \leq \frac{1}{q} \mu) \\
&\text{(Let $(1 - \delta) = \frac{1}{q}$  and apply (\ref{eq:cb:leq}))} \\
&\leq  \Big(\cfrac{e^{-\delta}}{(1 - \delta)^{(1 - \delta)}} \Big)^\mu \\
&=  \Big(e^{(\frac{1}{q} - 1)} q^{\frac{1}{q}} \Big)^{pk} = \Psi.
\end{align*}
 
\begin{bound}
\label{bd:chernoff}
By applying the Chernoff Bound, Q-error from random uniform sampling with replacement is bounded by:
\begin{align*}
\mathds{P}(\text{Q-error $\leq$ q}) 
&\geq \mathds{P}(\text{Q-error $<$ q})\\
&= 1 - \mathds{P}(\text{Q-error $\geq$ q})\\
&= 1 - \mathds{P}(Pred \geq q X) - \mathds{P}(Pred \leq \frac{1}{q} X) \\
&= 1-\Omega-\Psi\\
& \geq 1 -  \Big(\cfrac{e^{q - 1}}{q^{q}} \Big)^{pk} - \Big(e^{(\frac{1}{q} - 1)} q^{\frac{1}{q}} \Big)^{pk}.
\end{align*}
\end{bound}

\noindent\textbf{CE bounded by the Bernstein Inequality.}
The Bernstein inequality \cite[Reference][Proposition 1.4]{bardenet2015concentration} can be written as,
\begin{align}
\mathds{P}\Bigg(\cfrac{\sum\limits_{t=1}^k (x_t - \mu_x)}{k} > \epsilon \Bigg)  \leq \exp \Big(- \cfrac{k \epsilon^2}{ 2 \sigma^2 + 2 M \epsilon / 3} \Big),
\end{align}
where $M =\max\limits_{i=0}^n |x_i|$ and $\epsilon > 0$. In our Bernoulli case, $M = 1$, and $\mu_x = p$ is the population mean. Therefore, we have

\begin{align*}
\mathds{P}(Pred > q X) &= \mathds{P}(\hat{X} > kpq) \\
    &= \mathds{P}(\hat{X} - kp > kpq - kp) \\
    &= \mathds{P}\Bigg( \frac{\sum\limits_{t=1}^k (X_t - p)}{k} > pq - p\Bigg) \\
    &\text{(Let $\epsilon = pq - p$)} \\
    &= \mathds{P}\Bigg( \frac{\sum\limits_{t=1}^k (X_t - p)}{k} > \epsilon\Bigg) \\
    &\leq \exp \Big(- \cfrac{k \epsilon^2}{ 2 \sigma^2 + 2 M \epsilon / 3} \Big) =   \exp \Big(- \frac{k (pq-p)^2}{2 \sigma^2  + 2 (pq-p) / 3}\Big)= \Omega,
\end{align*}

\begin{align*}
    \mathds{P}(Pred < X/q)
    &= \mathds{P}(\hat{X} < \frac{kp}{q}) \\
    &= \mathds{P}(k - \hat{Y} < \frac{kp}{q}) \\
    &= \mathds{P}(\hat{Y} > k - \frac{kp}{q}) \\
    &= \mathds{P}(\hat{Y} - k\bar{p} > k - \frac{kp}{q} - k\bar{p})  \\
    &= \mathds{P}\Bigg( \frac{\sum\limits_{t=1}^k (Y_t - \bar{p})}{k} > 1- p/q - (1-p)\Bigg) \\
    &= \mathds{P}\Bigg( \frac{\sum\limits_{t=1}^k (Y_t - \bar{p})}{k} > p - p/q \Bigg) \\
    &\text{(Let $\epsilon = p - p/q$)} \\
    &= \mathds{P}\Bigg( \frac{\sum\limits_{t=1}^k (Y_t - \bar{p})}{k} > \epsilon \Bigg) \\
    &\leq \exp \Big(- \cfrac{k \epsilon^2}{ 2 \sigma^2 + 2 M \epsilon / 3} \Big) =  \exp \Big(- \frac{k (p-p/q)^2}{2 \sigma^2 + 2 (p-p/q) / 3} \Big)= \Psi.
\end{align*}
 
\begin{bound} 
\label{bd:bernstein}
By applying the Bernstein's Inequality, the Q-error from random uniform sampling with replacement is bounded by:
 
$\mathds{P}(\text{Q-error $\leq$ q})
 \geq 1-\Omega-\Psi = 1 -   \exp \Big(- \frac{k (pq-p)^2}{2 \sigma^2 + 2 (pq-p) / 3} \Big)  -  \exp \Big(- \frac{k (p-p/q)^2}{2 \sigma^2 + 2 (p-p/q) / 3} \Big)$.
\end{bound}
 
\noindent Putting together the $\Omega$ and $\Psi$ in the above bounds, we can derive Theorem \ref{theory:with_replacement}. 
 
\begin{figure}[h!]
    \centering
     \begin{subfigure}[b]{0.32\textwidth}
        \includegraphics[width=1\textwidth]{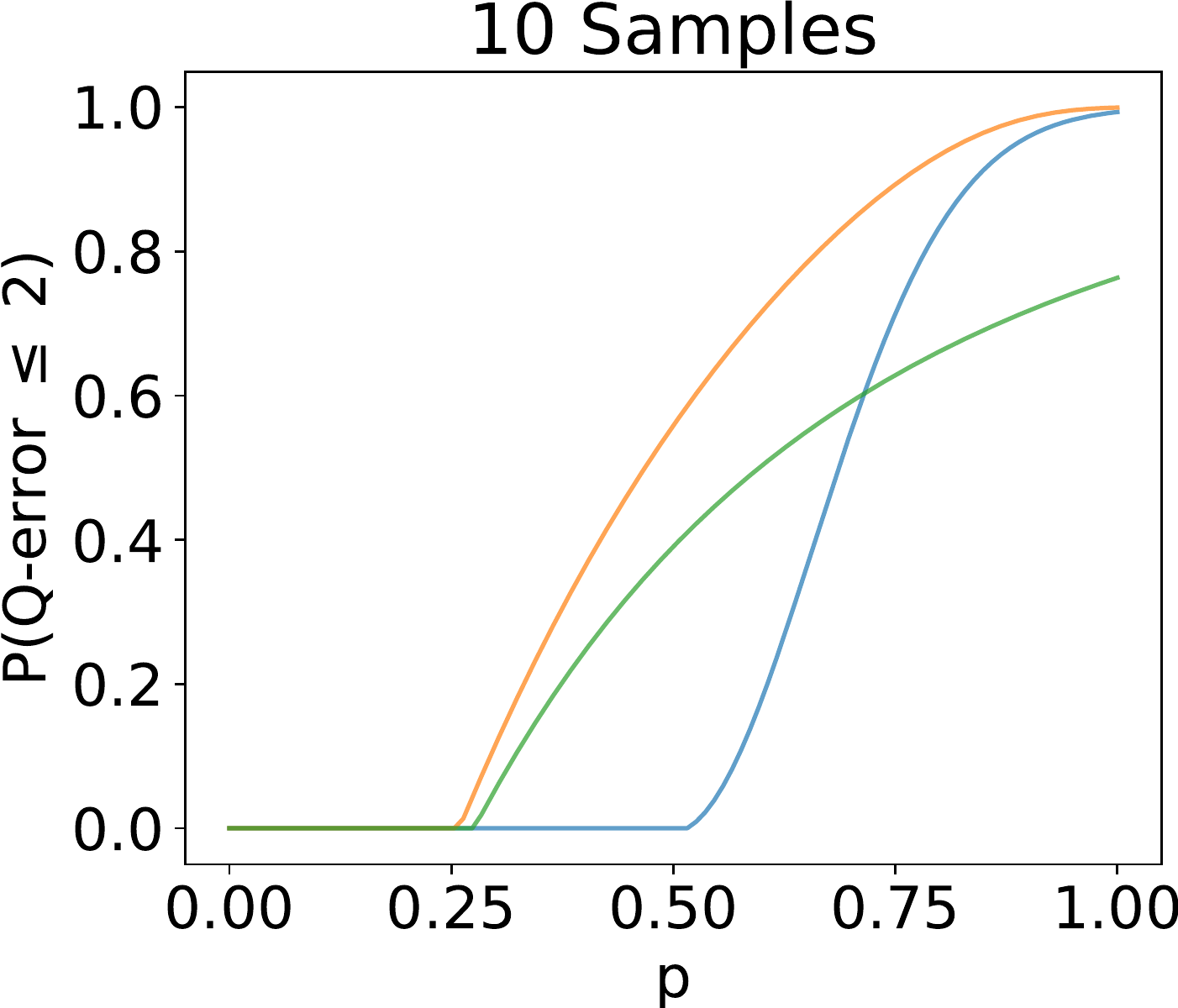}
    \end{subfigure}
     \begin{subfigure}[b]{0.32\textwidth}
        \includegraphics[width=1\textwidth]{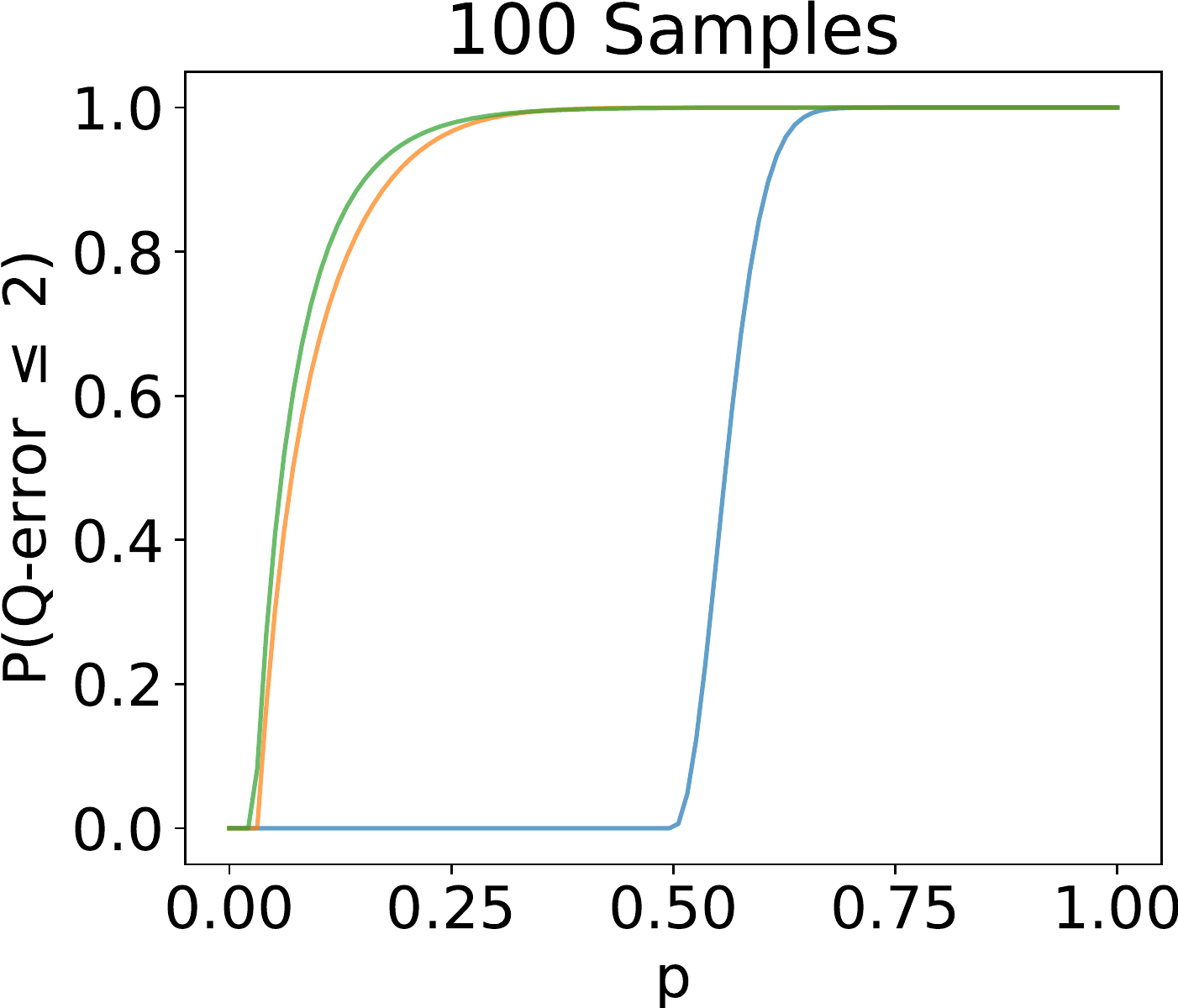}
    \end{subfigure}
    \begin{subfigure}[b]{0.32\textwidth}
        \includegraphics[width=1\textwidth]{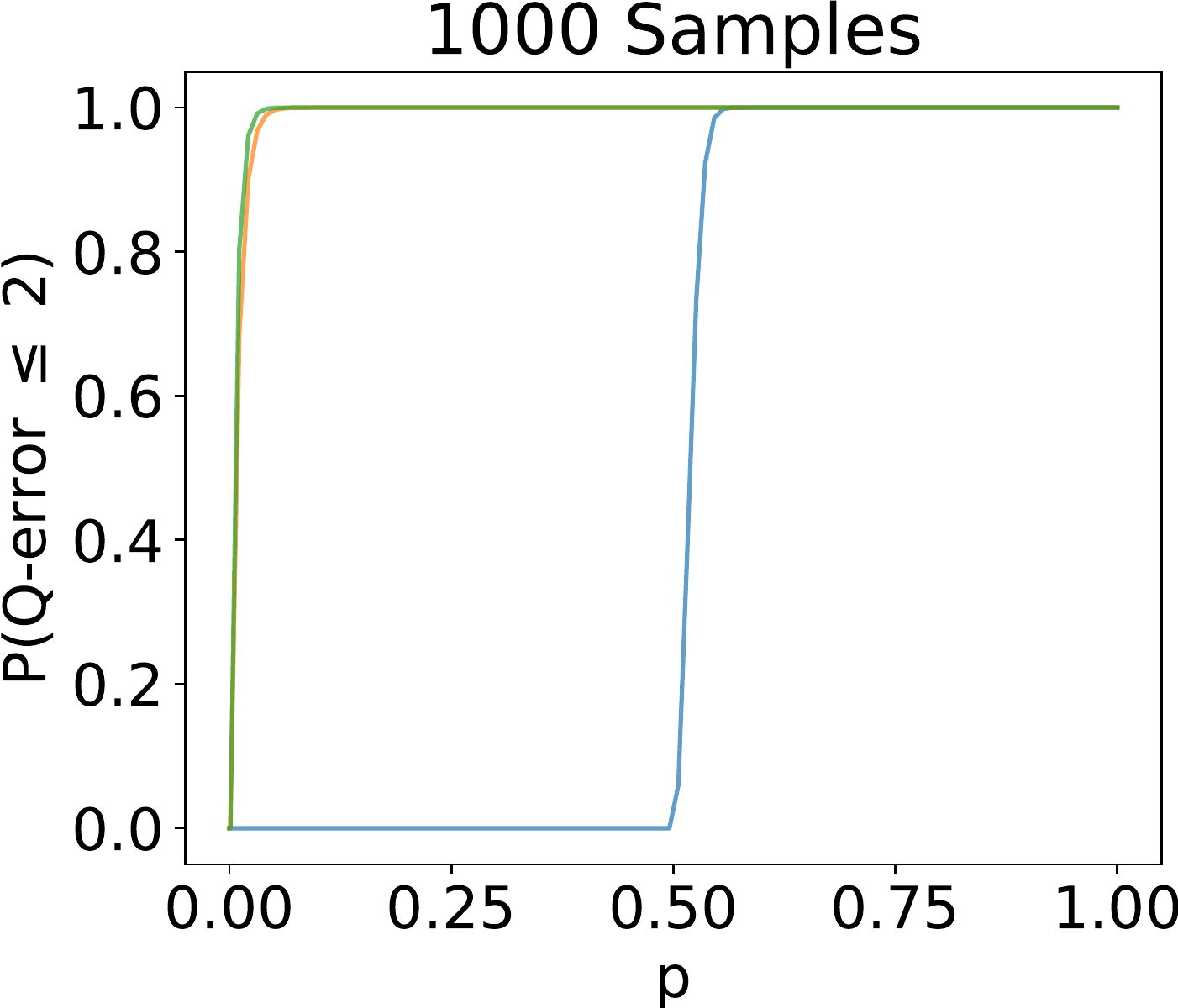}
    \end{subfigure} \vspace{-0.1in}
    \begin{subfigure}[b]{0.70\textwidth}
        \includegraphics[width=1\textwidth]{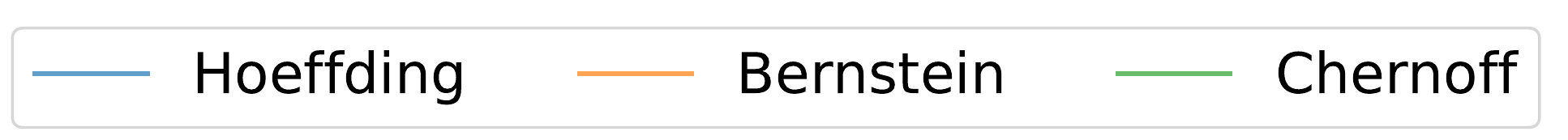}
    \end{subfigure}
    \caption{Plotting and comparing the Chernoff, Bernstein's Inequality and Hoeffding's Inequality bounds (see Appendix). }
    \label{fig:chernoff_bernstein_hoeffding}
\end{figure}

\vspace{0.1in}\noindent\textbf{Visualization.}
Figure \ref{fig:chernoff_bernstein_hoeffding} plots and compares the bounds shown in this section. With only 100 samples (rows), Bound \ref{bd:chernoff} and \ref{bd:bernstein} already demonstrate some tightness: the Q-error is likely to be small (more than 80\% chance with $q\le2$) when a query predicate has a cardinality of $p=0.2$; with 1K samples, the Q-error is almost always small ($q\le$2). Figures 2 and 3 further demonstrate the probability with different $p$ values and the Q-error with 95\% confidence. 
 
We also show in Figure \ref{fig:3d_replace} the 3-D plotting of the probability that random sampling has a Q-error that is better than a given  threshold. At a small sample size such as 100 or 1K rows irrelevant to the table size, random uniform sampling with replacement already provides a promising q-error and has a tight bound.

\begin{figure}[h!]
    \centering
     \begin{subfigure}[b]{0.42\textwidth}
        \includegraphics[width=1\textwidth]{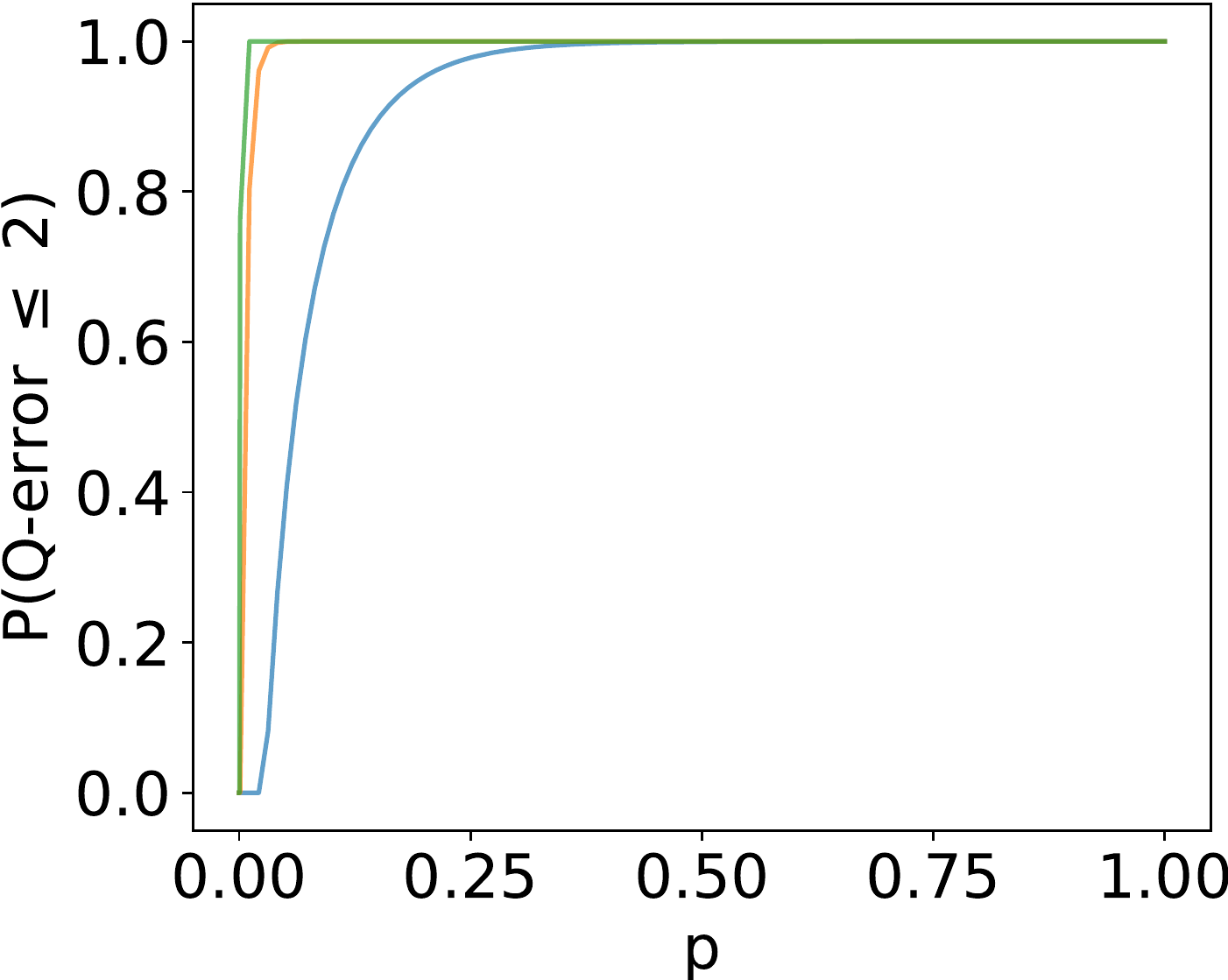}
    \end{subfigure}
    \begin{subfigure}[b]{0.42\textwidth}
        \includegraphics[width=1\textwidth]{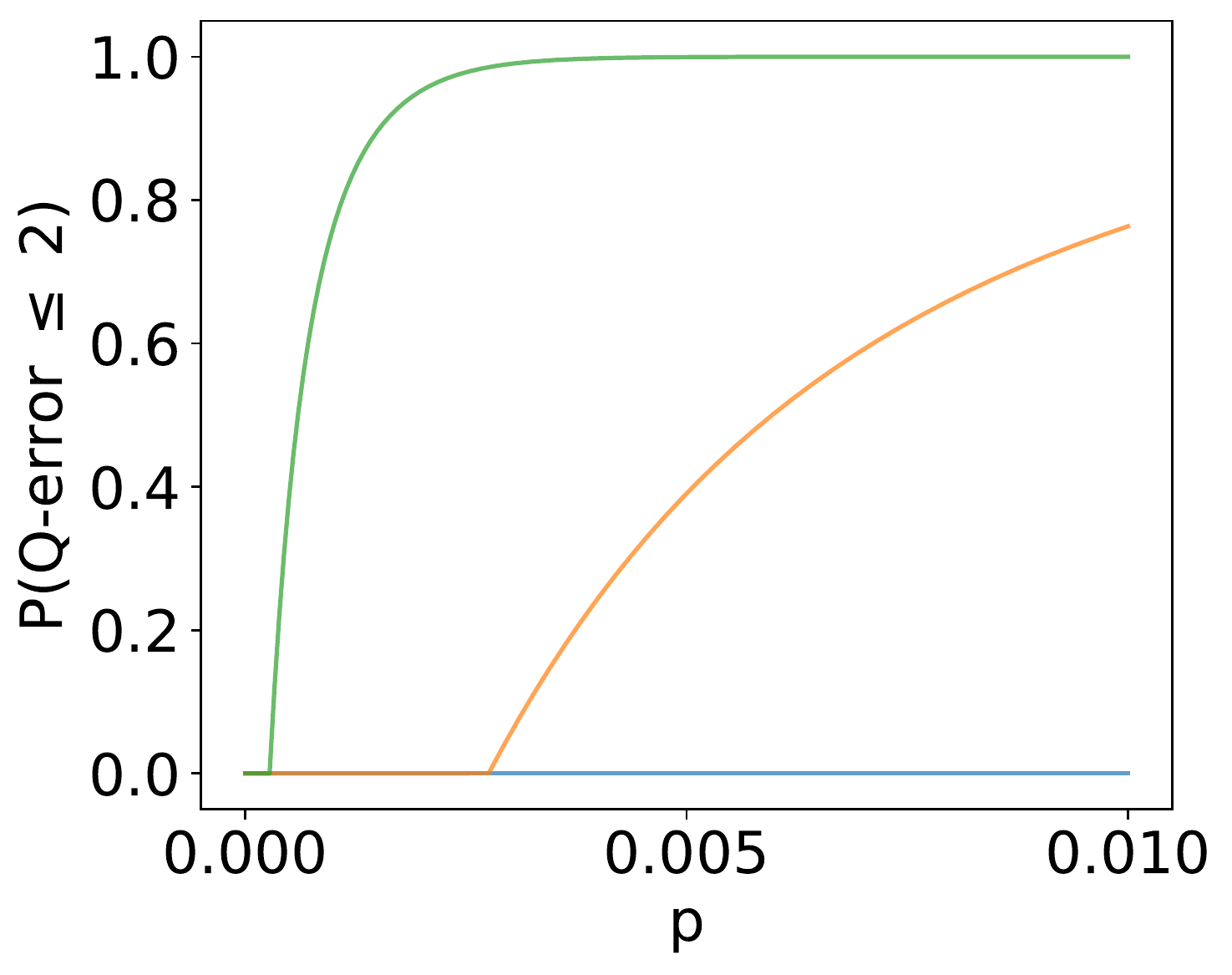}
    \end{subfigure}
    \begin{subfigure}[b]{0.70\textwidth}
        \includegraphics[width=1\textwidth]{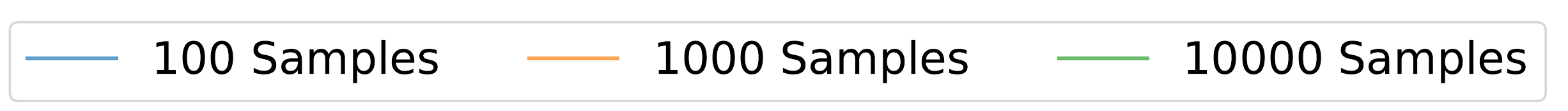}
    \end{subfigure} \vspace{-0.1in}
    \caption{Plotting the bounds in terms of $p$.
     (Left) Results for $p \in (0, 100\%]$. (Right) Results for $p\in (0, 1\%]$. }
     \label{fig:n_samples_with_replacement}
\end{figure}

\begin{figure}[h!]
    \centering
     \begin{subfigure}[b]{0.42\textwidth}
        \includegraphics[width=1\textwidth]{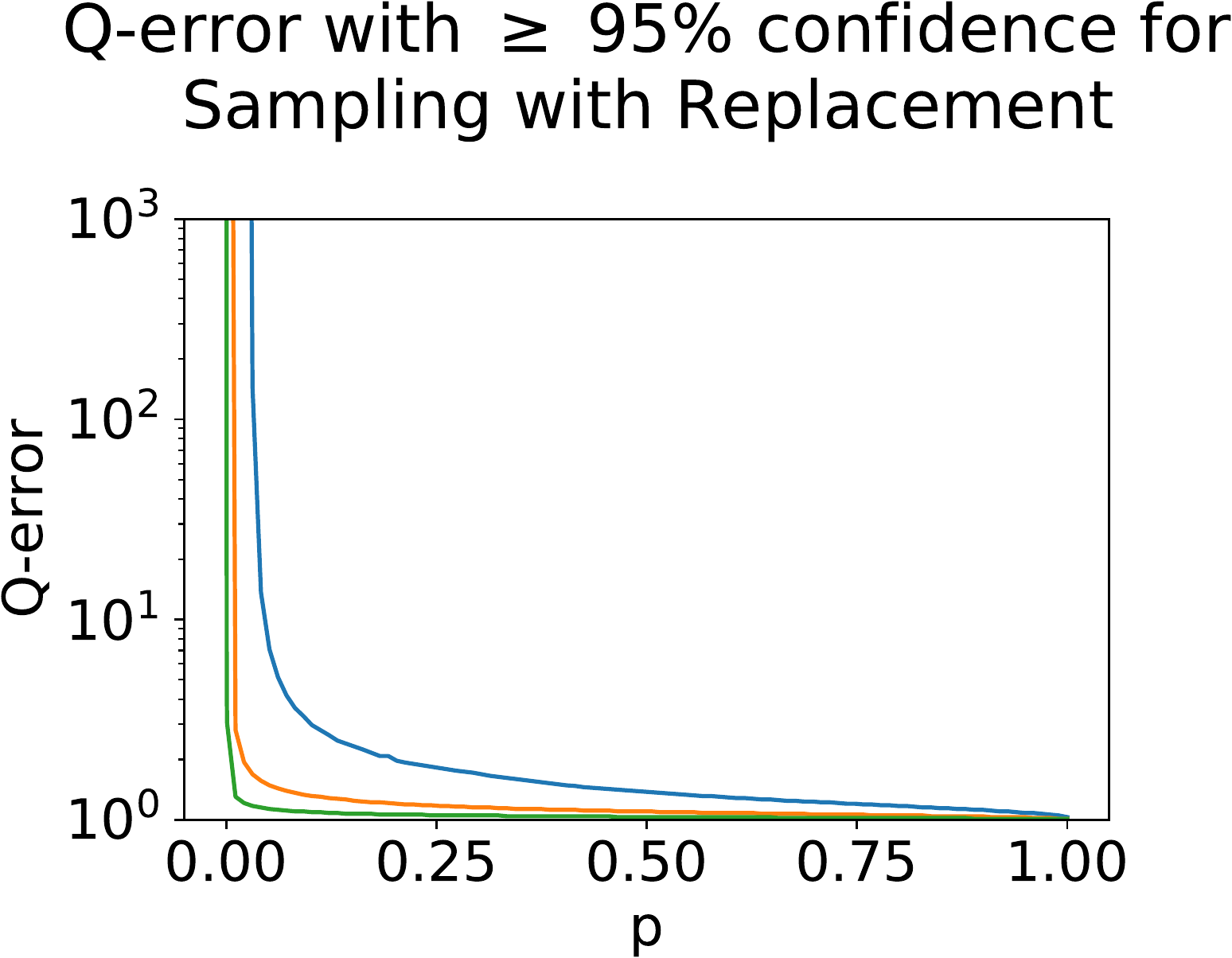}
    \end{subfigure}
     \begin{subfigure}[b]{0.42\textwidth}
        \includegraphics[width=1\textwidth]{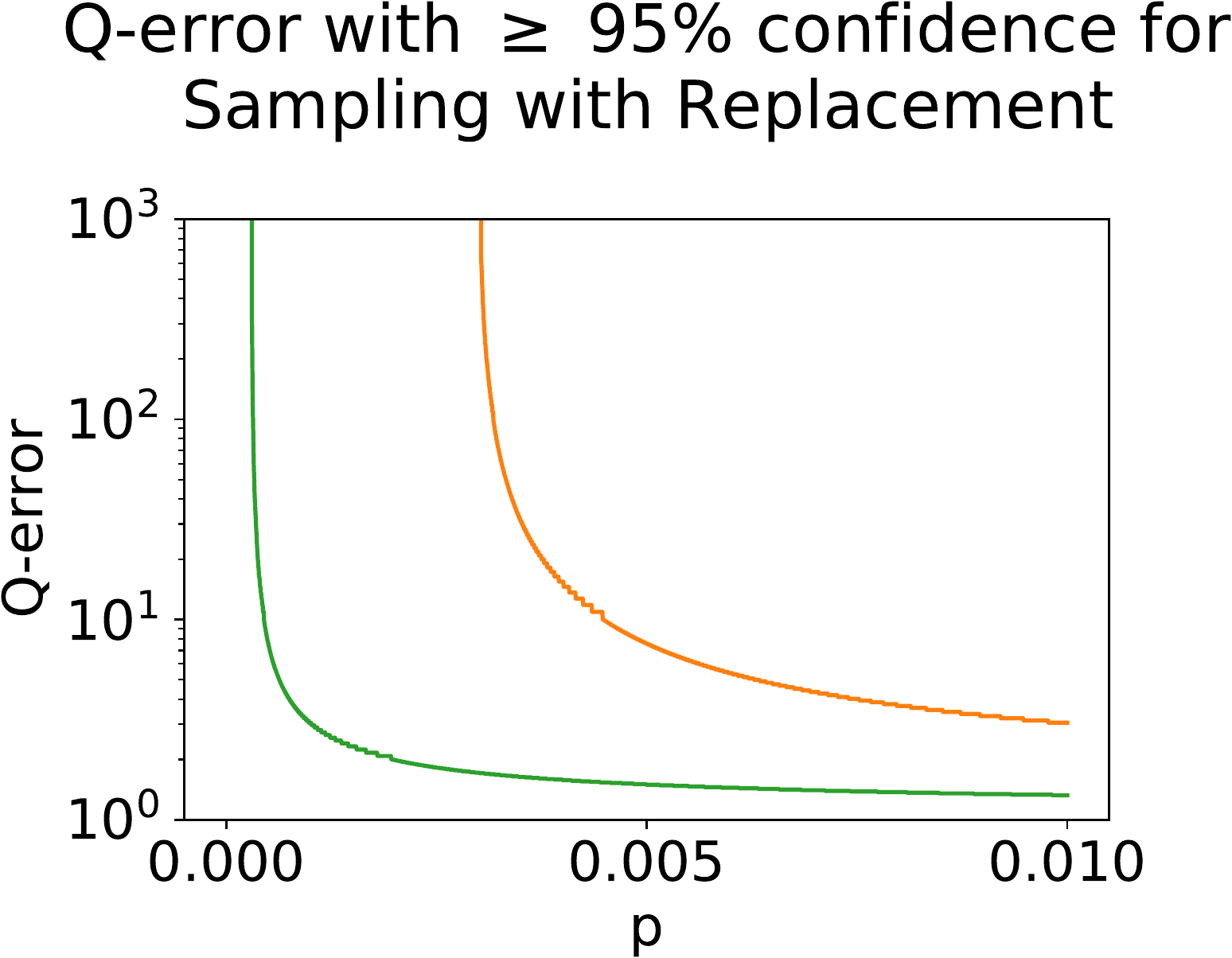}
    \end{subfigure}
    \begin{subfigure}[b]{0.70\textwidth}
        \includegraphics[width=1\textwidth]{fig/legend_n_samples.pdf}
    \end{subfigure} \vspace{-0.1in}
    \caption{Plotting the Q-error with 95\% confidence in terms of $p$.
     (Left) Results for $p \in (0, 100\%]$. (Right) Results for $p\in (0, 1\%]$.
    }
    \label{fig:n_samples_with_replacement_vs_q}
\end{figure}

\begin{figure*}[h!]
    \centering
     \begin{subfigure}[b]{0.32\textwidth}
        \includegraphics[width=1\textwidth]{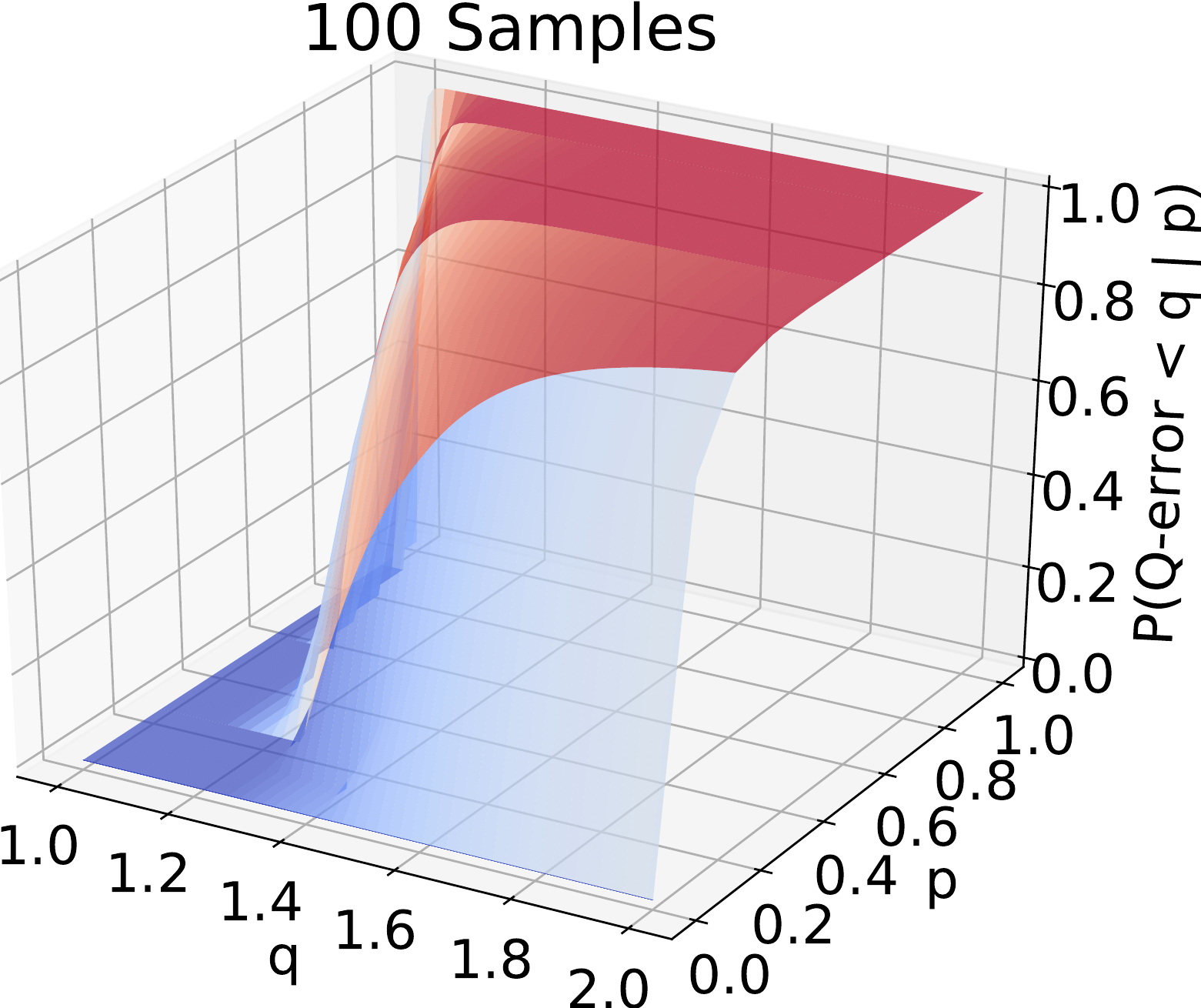}
    \end{subfigure}
    \begin{subfigure}[b]{0.32\textwidth}
        \includegraphics[width=1\textwidth]{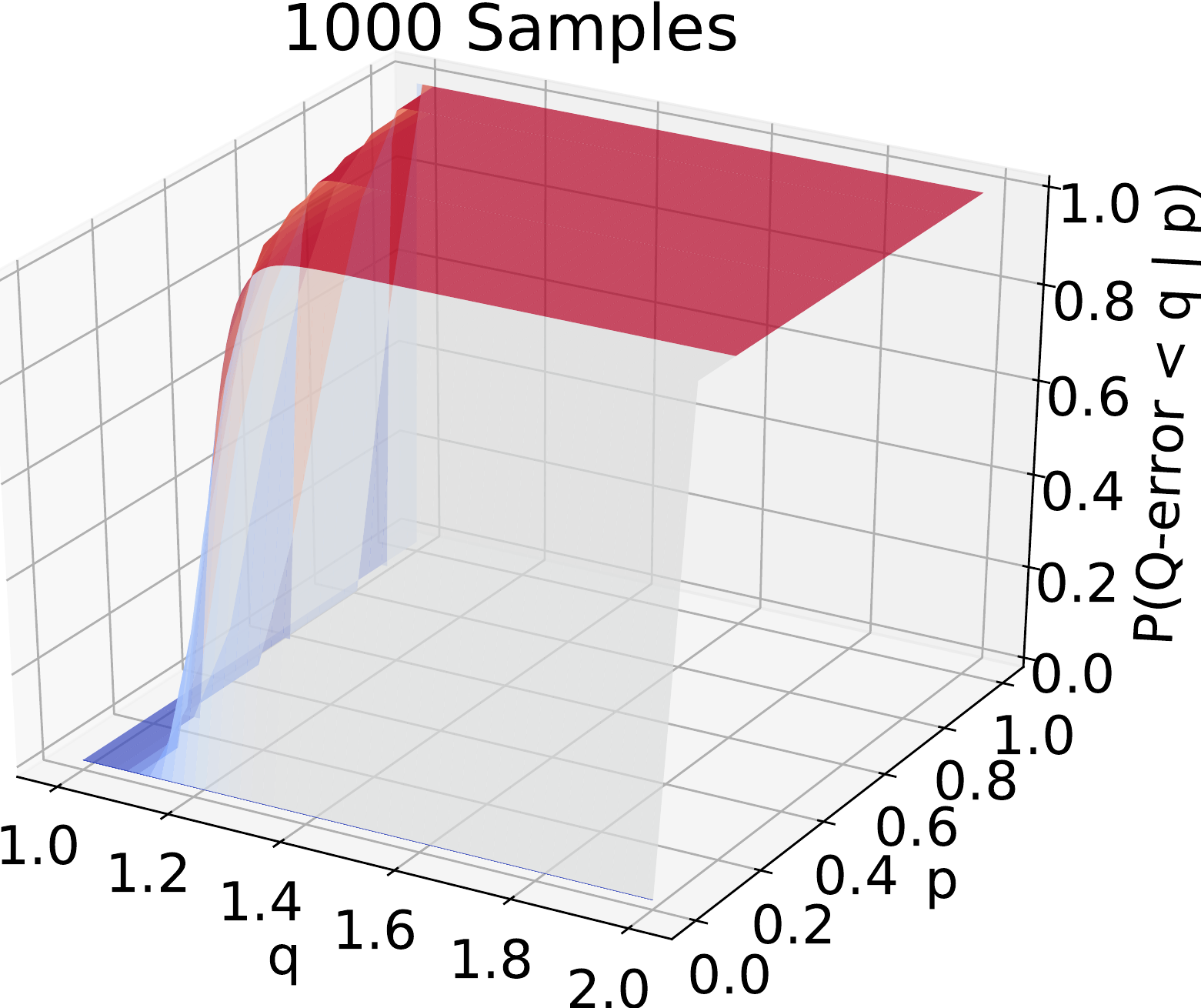}
    \end{subfigure}
    \begin{subfigure}[b]{0.32\textwidth}
        \includegraphics[width=1\textwidth]{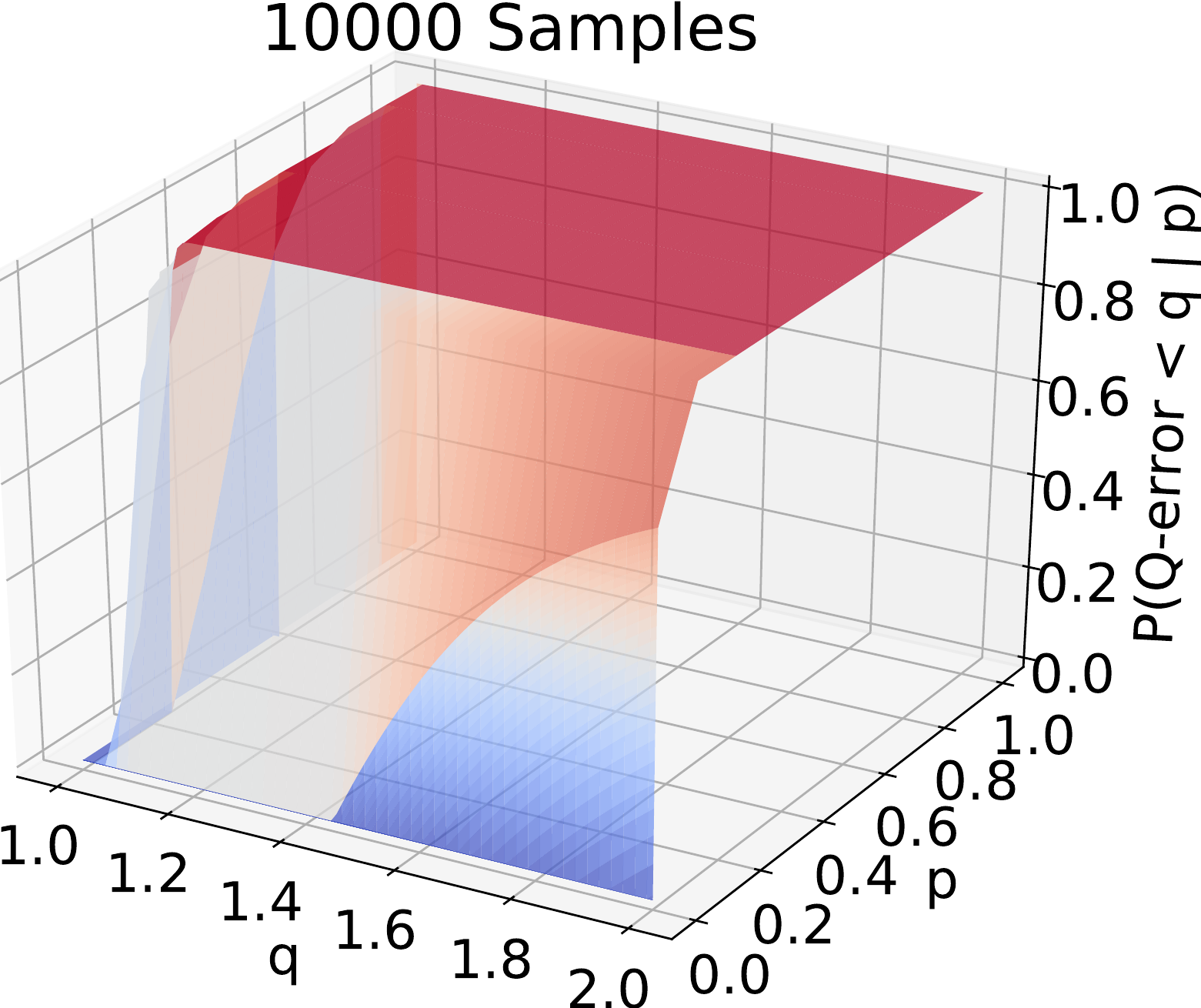}
    \end{subfigure} 
    \caption{3D plotting of the probability that random uniform sampling's error (z-axis) is better than $q$ with a predicate cardinality $p$. }
    \label{fig:3d_replace}
\end{figure*}

\section{Bounds of Sampling without Replacement}
\label{sec:rs_np}
 
When samples are drawn uniformly at random without replacement, we assume there are at least 2 rows.  The Hoeffding-Serfling inequality and the Bernstein-Serfling inequality developed in a recent study \citep{bardenet2015concentration} founded Theorem \ref{theory:without_replacement} for random uniform sampling without replacement. Adding Serfling's results further gives a tighter bound. 
 
\begin{theorem}
\label{theory:without_replacement}
For cardinality estimation over single tables, the Q-error of random uniform sampling without replacement is bounded by
 
\begin{align*}
& \mathds{P}(\text{Q-error} \leq q)  \geq  1 - \Omega -\Psi \text{, where}\\
\Omega &= \min \Bigg(2 \exp \Big(- \frac{k}{\zeta^2} (- \sqrt{ 2 \zeta \rho \sigma^2 (pq - p)  + \rho^2 \sigma^4} + (pq - p)\zeta + \sigma^2 \rho ) \Big),  \\
&~~~~~~~~~~~~~~ \exp \Big(-\cfrac{2 k (pq - p)^2}{\rho} \Big) \Bigg),\\
\Psi &= \min \Bigg(2 \exp \Big(- \frac{k}{\zeta^2} (- \sqrt{ 2 \zeta \rho \sigma^2 (p - p/q)  + \rho^2 \sigma^4} + (p - p/q)\zeta + \sigma^2 \rho  ) \Big), \\
&~~~~~~~~~~~~~~ \exp \Big(- \cfrac{2k(p - p/q)^2}{\rho} \Big) \Bigg).
\end{align*}
 
\end{theorem}

\noindent\textbf{CE bounded by Hoeffding-Serfling.} For convenience, we define
\begin{center}
$\rho = \begin{cases}
    1 - (k - 1)/n & \text{if $k \leq n / 2$} \\ 
    (1 - k/n) ( 1 + 1/k) & \text{if $k > n / 2$}
    \end{cases}
    $
, and
$\zeta = \begin{cases}
    4/3 +  \sqrt{\frac{k (k-1)}{n (n-k+1)}} & \text{if $k \leq n / 2$} \\ 
    4/3 + \sqrt{\frac{(n-k-1)(n-k)}{(k+1)n}} & \text{if $k > n / 2$}.
    \end{cases}
$
\end{center} Let $\chi = (r_1,...,r_n)$ be a finite population of $n$ rows. We sample $k$ rows (i.e. $x_1, ..., x_k$) without replacement from the population such that $k < n$. Denote $\mathds{E}[r_i] =  \mu_x$, $a = \min\limits_{1 \leq i \leq n} r_i$, and $b = \max\limits_{1 \leq i \leq n} r_i$.
The Hoeffding-Serfling inequality \cite[Reference][Corollary 2.5]{bardenet2015concentration}  can be written as:
 
\begin{center}
$\mathds{P}\Bigg(\cfrac{\sum\limits_{t=1}^k (x_t - \mu_x)}{k} > (b - a) \sqrt{\cfrac{\rho \log(1/\delta)}{2k}}\Bigg)  \leq \delta $,
\end{center}

\noindent where $\delta \in [0, 1]$. For Bernoulli distribution, we can simplify these two inequalities by $a = 0$, $b = 1$, and $\sigma^2 = p (1-p)$. 
 
Let $\delta =  \exp \Big(- \frac{2k(pq - p)^2}{\rho} \Big)$. 
Denote $Y_i = \neg X_i$ such that $Y_i = 0$ if the sampled row satisfy the predicate and $y_i = 1$ if the row does not satisfy the predicate. We have $\hat{Y} = \sum\limits_{i=0}^k\hat y_i$ on the sample. It is easy to see $\hat{X} = k - \hat{Y}$. We denote $\bar{p} = \mathds{E}[Y_i] = 1 - p$. The variances for $Y$ is $\bar{p} (1-\bar{p}) = (1-p)p$. By applying the Hoeffding-Serfling inequality, we have the over-estimation probability as
 
\begin{align*}
\mathds{P}(Pred > q X) &= \mathds{P}(\hat{X} > kpq) \\
    &= \mathds{P}(\hat{X} - kp > kpq - kp) \\
    &= \mathds{P}\Bigg( \frac{\sum\limits_{t=1}^k (X_t - p)}{k} > pq - p\Bigg) \\
    &= \mathds{P}\Bigg( \frac{\sum\limits_{t=1}^k (X_t - p)}{k} > \sqrt{\cfrac{\rho \log(1/\delta)}{2k}}\Bigg) \\
    &\leq \delta  =  \exp \Big(- \frac{2k(pq - p)^2}{\rho} \Big)=\Omega.
\end{align*}
 
\noindent Let $\delta = \exp \Big(- \frac{2k(p - p/q)^2}{\rho} \Big)$. We have the under-estimation probability as
 
\begin{align*}
    \mathds{P}(Pred < X/q)
    &= \mathds{P}(\hat{X} < \frac{kp}{q}) \\
    &= \mathds{P}(k - \hat{Y} < \frac{kp}{q}) \\
    &= \mathds{P}(\hat{Y} > k - \frac{kp}{q}) \\
    &= \mathds{P}(\hat{Y} - k\bar{p} > k - \frac{kp}{q} - k\bar{p})  \\
    &= \mathds{P}\Bigg( \frac{\sum\limits_{t=1}^k (Y_t - \bar{p})}{k} > 1- p/q - (1-p)\Bigg) \\
    &= \mathds{P}\Bigg( \frac{\sum\limits_{t=1}^k (Y_t - \bar{p})}{k} > p - p/q \Bigg) \\
    &\leq \delta = \exp \Big(- \frac{2k(p - p/q)^2}{\rho} \Big)=\Psi.
\end{align*}

\begin{bound} 
\label{bd:hs}
By applying the Hoeffding-Serfling Inequality, the Q-error of random uniform sampling without replacement is bounded by
\begin{align*}
\mathds{P}(\text{Q-error $\leq$ q}) & \geq  1-\Omega-\Psi = 1  - \exp \Big(-\cfrac{2 k (pq - p)^2}{\rho} \Big) - \exp \Big(- \cfrac{2k(p - p/q)^2}{\rho} \Big).\\ 
\end{align*}
\end{bound}
 
\noindent\textbf{CE bounded by the Bernstein-Serfling Inequality}.
Bernstein-Serfling inequality is usually tighter than Hoeffding-Serfling inequality unless $p$ is large. The Bernstein-Serfling inequality \cite[Reference][Corollary 3.6]{bardenet2015concentration}  can be written as:
 
\begin{center}
$\mathds{P}\Bigg(\cfrac{\sum\limits_{t=1}^k (x_t - \mu_x)}{k} > \sigma \sqrt{\cfrac{2 \rho \log(1/\delta)}{k}} + \cfrac{\zeta  (b-a) \log (1/\delta)}{k}\Bigg)  \leq  2\delta $,
\end{center}
 
\begin{figure}[h!]
    \centering
     \begin{subfigure}[b]{0.32\textwidth}
        \includegraphics[width=1\textwidth]{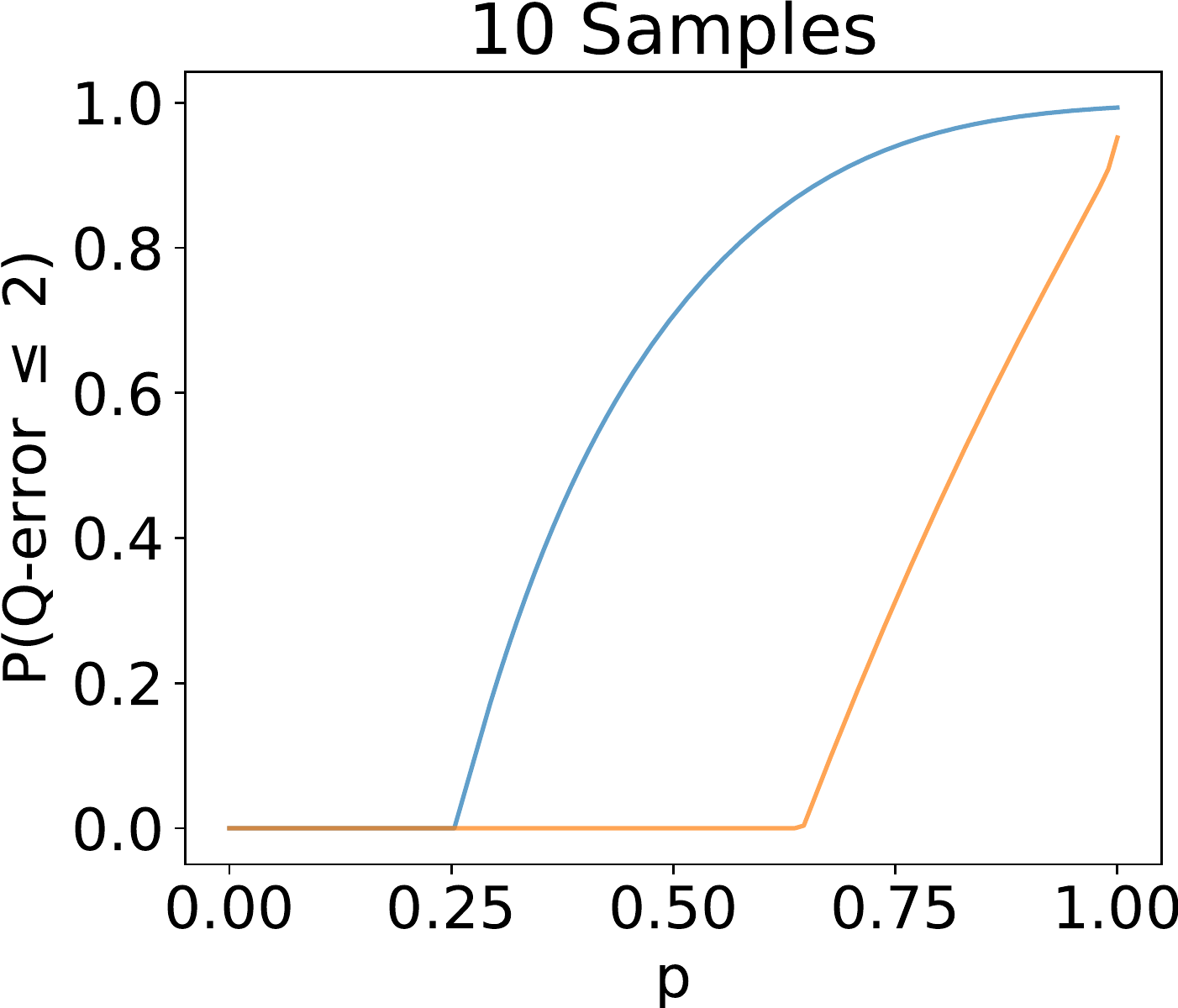}
    \end{subfigure}
     \begin{subfigure}[b]{0.32\textwidth}
        \includegraphics[width=1\textwidth]{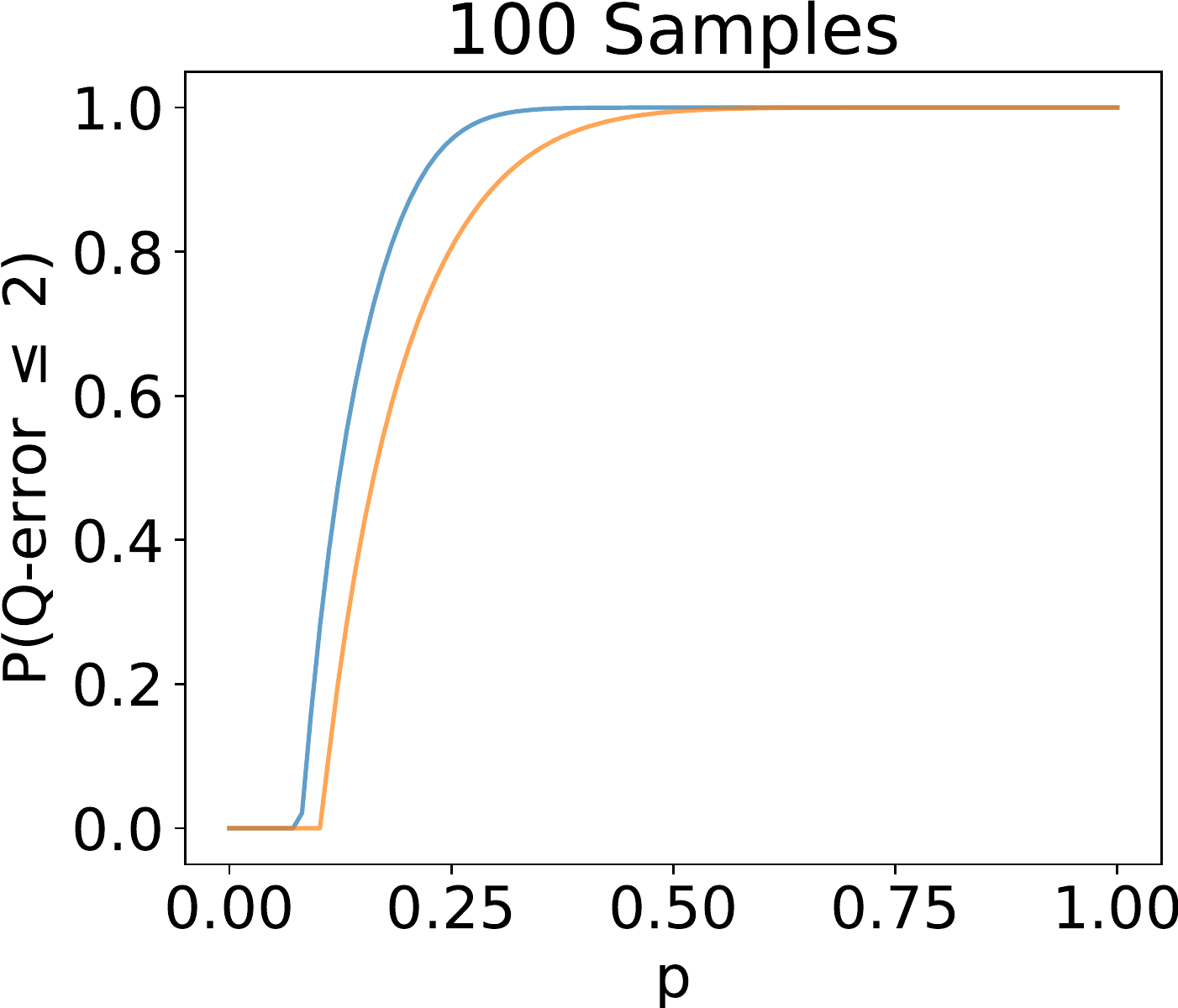}
    \end{subfigure}
    \begin{subfigure}[b]{0.32\textwidth}
        \includegraphics[width=1\textwidth]{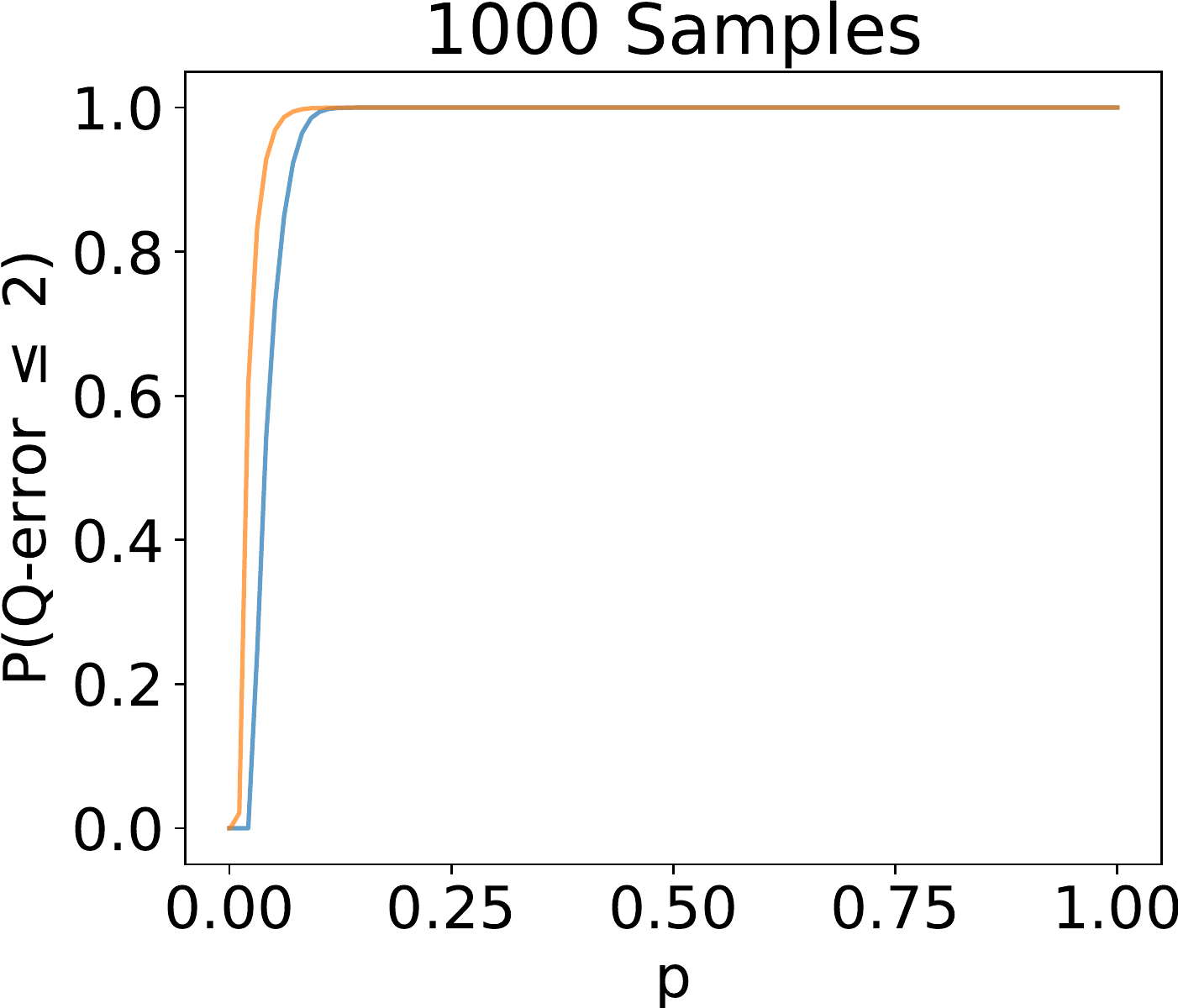}
    \end{subfigure} \vspace{-0.1in}
    \begin{subfigure}[b]{0.70\textwidth}
        \includegraphics[width=1\textwidth]{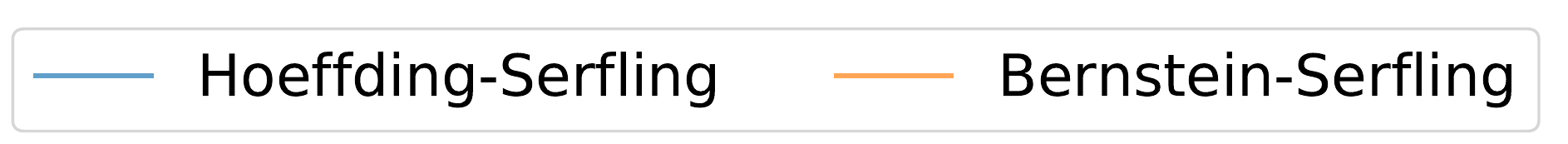}
    \end{subfigure}
    \caption{Comparing the Hoeffding-Serfling and Bernstein-Serfling bounds, assuming the length of the table is 1 billion. }
    \label{fig:hoeffding_bernstein_serfling}
\end{figure}

\begin{figure*}[htb]
    \centering
     \begin{subfigure}[b]{0.32\textwidth}
        \includegraphics[width=1.01\textwidth]{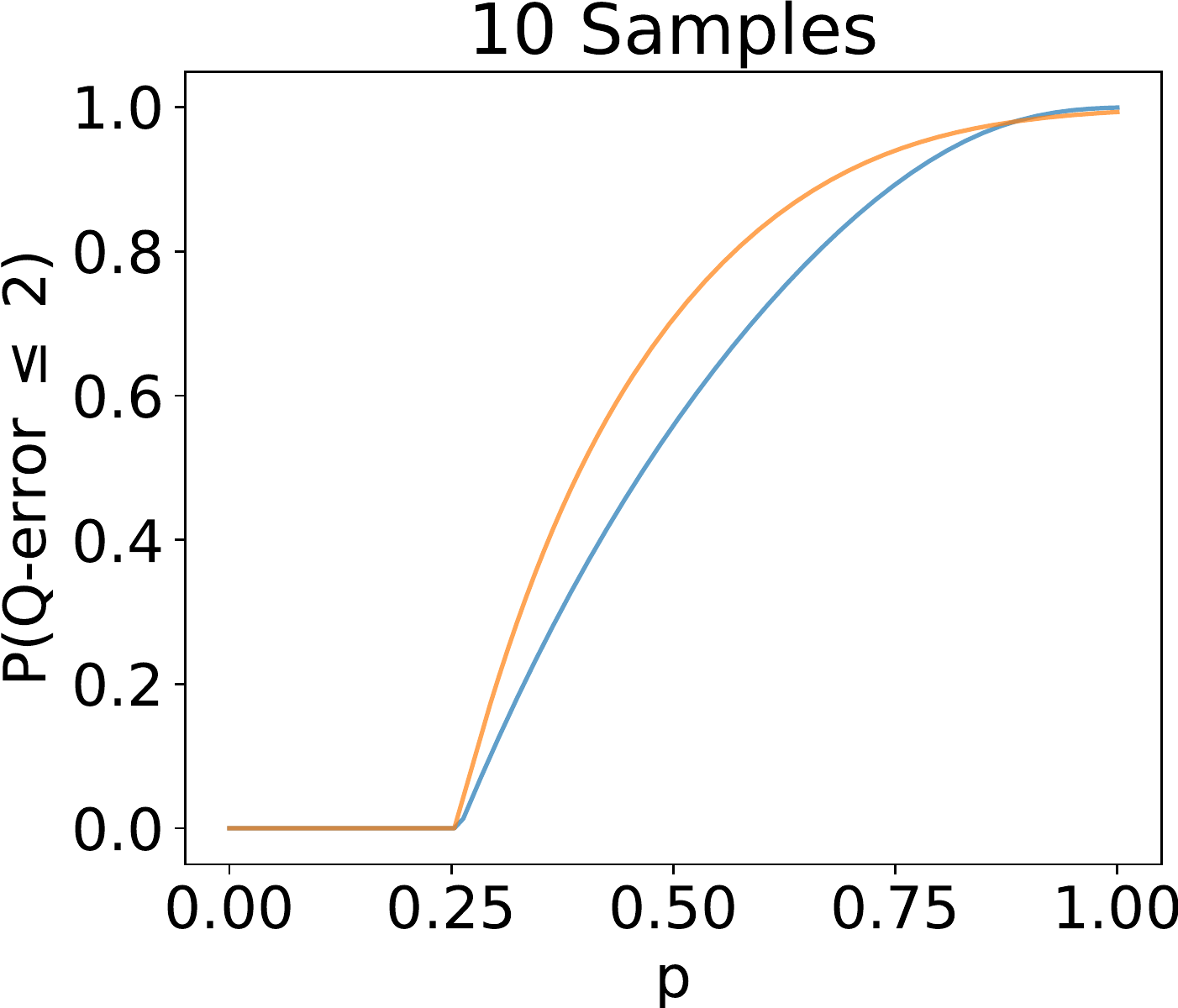}
    \end{subfigure}
    \begin{subfigure}[b]{0.32\textwidth}
        \includegraphics[width=1.01\textwidth]{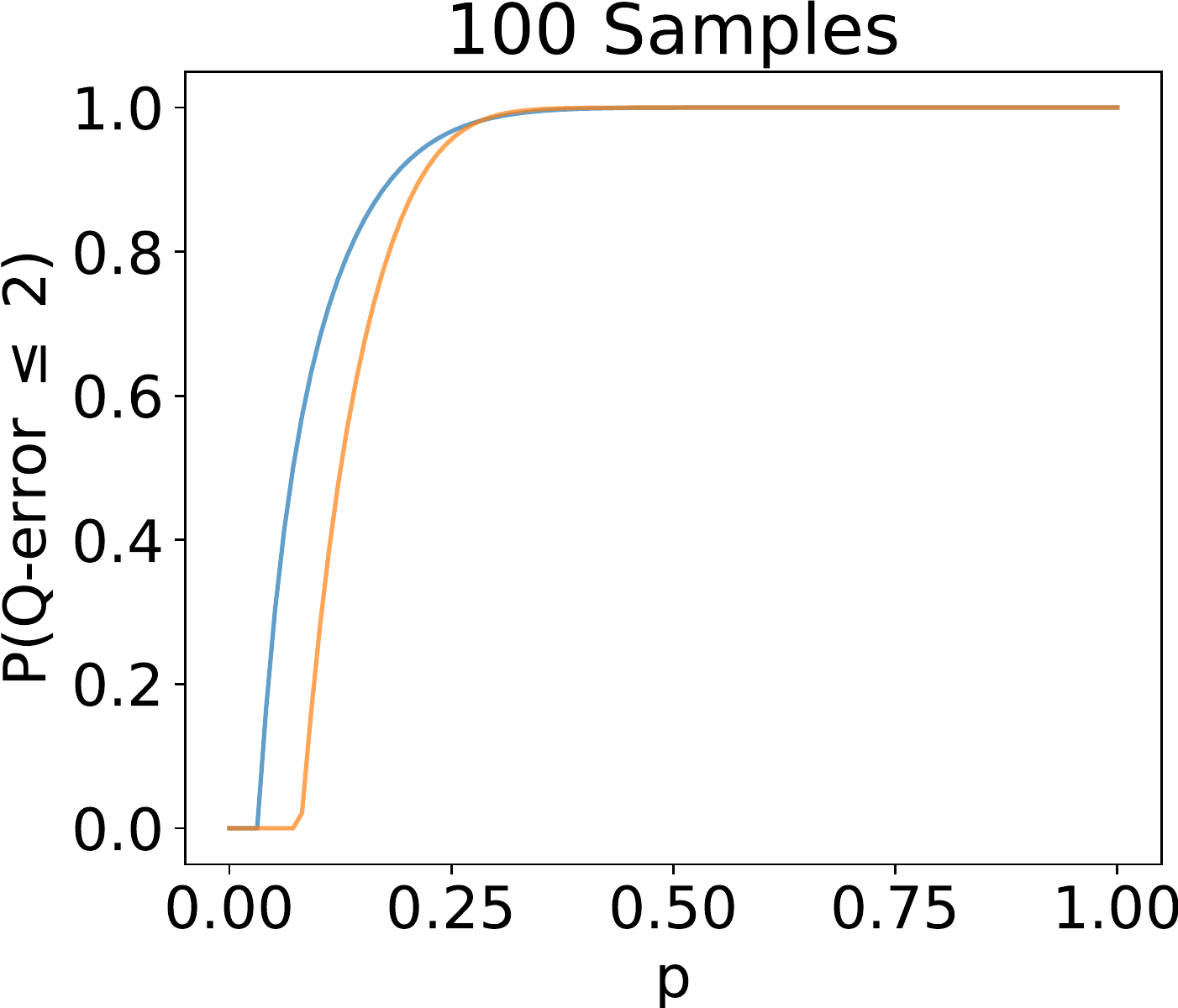}
    \end{subfigure}
    \begin{subfigure}[b]{0.32\textwidth}
        \includegraphics[width=1.01\textwidth]{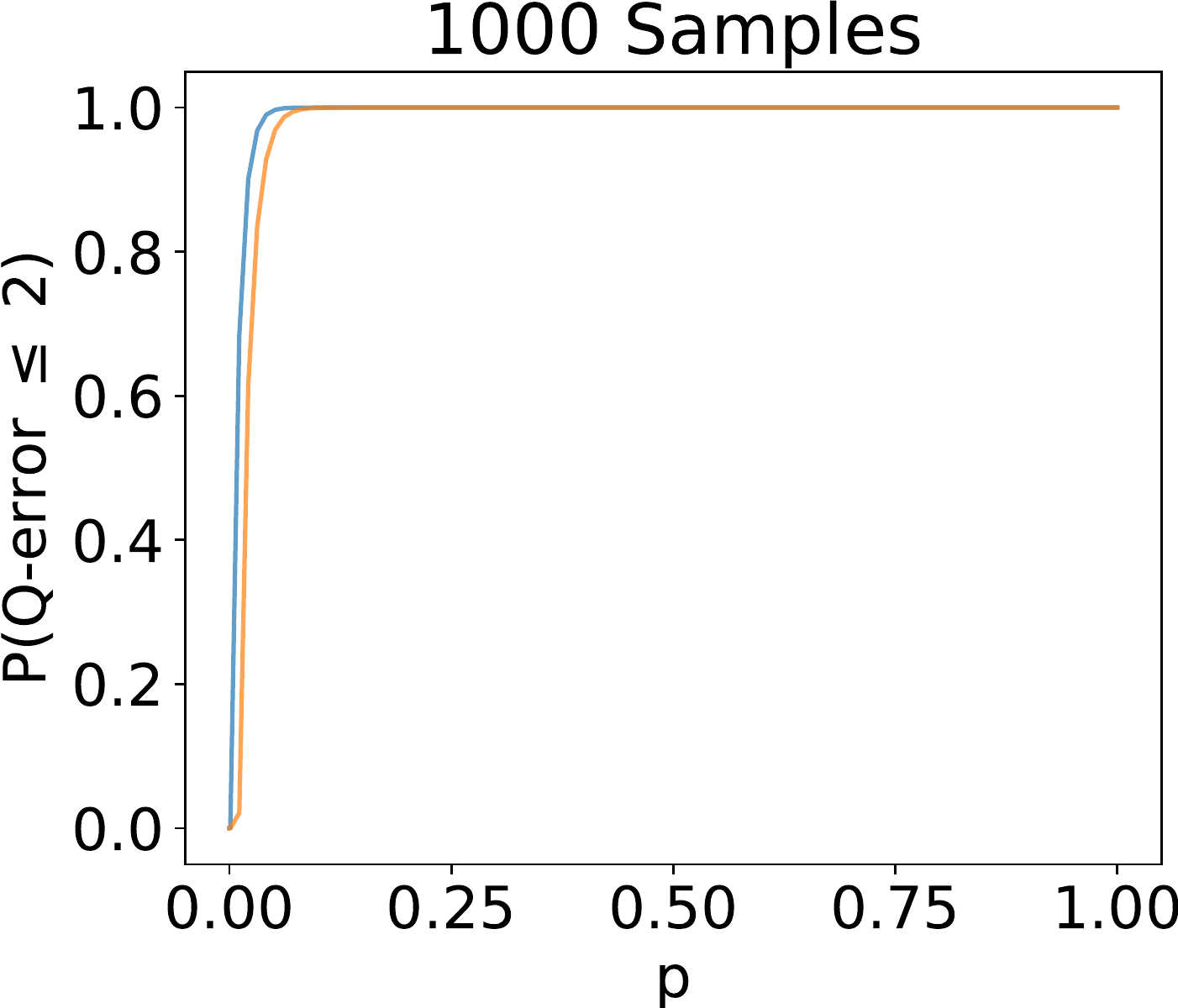}
    \end{subfigure}
     \begin{subfigure}[b]{0.32\textwidth}
        \includegraphics[width=1.01\textwidth]{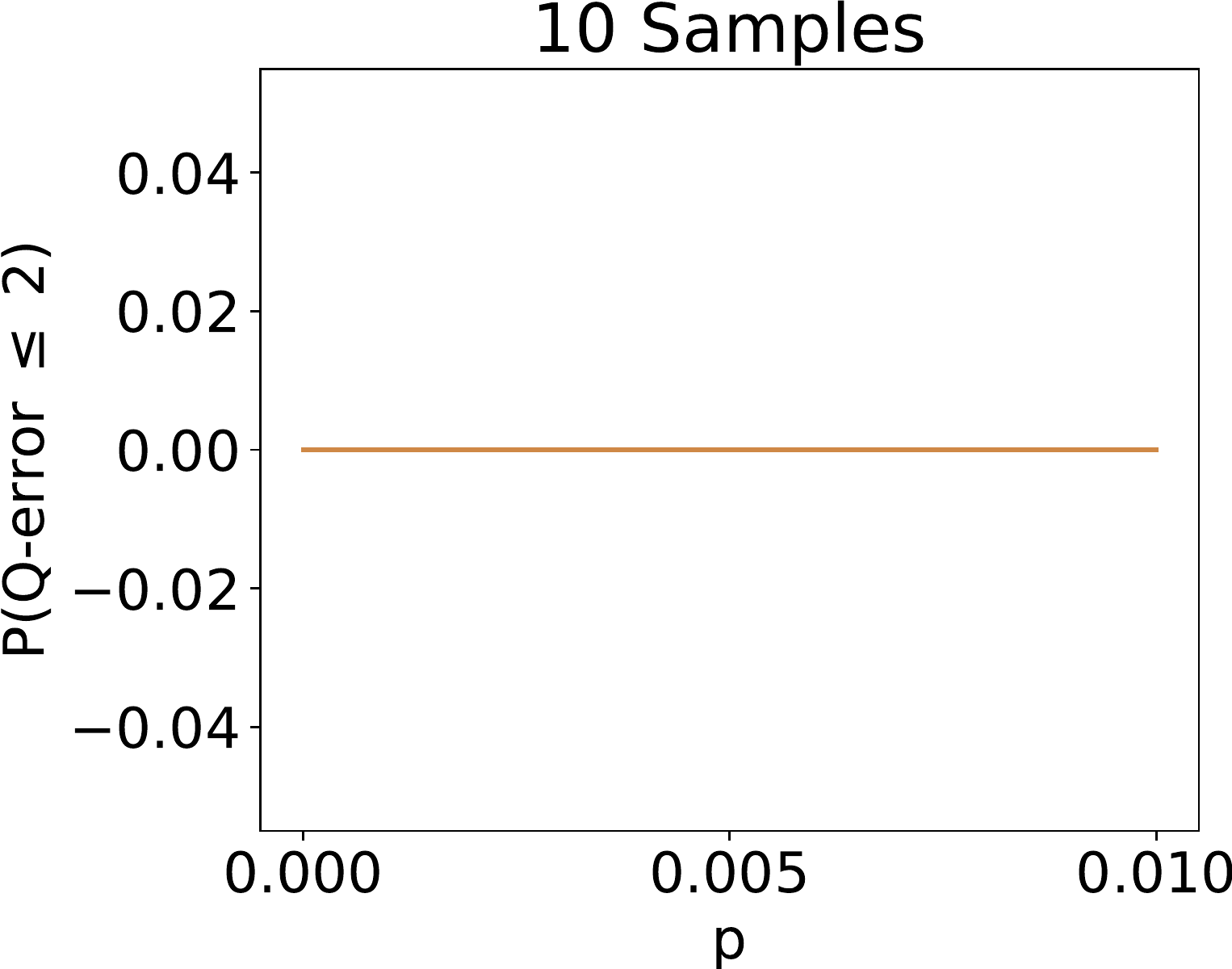}
    \end{subfigure}
    \begin{subfigure}[b]{0.32\textwidth}
        \includegraphics[width=1.01\textwidth]{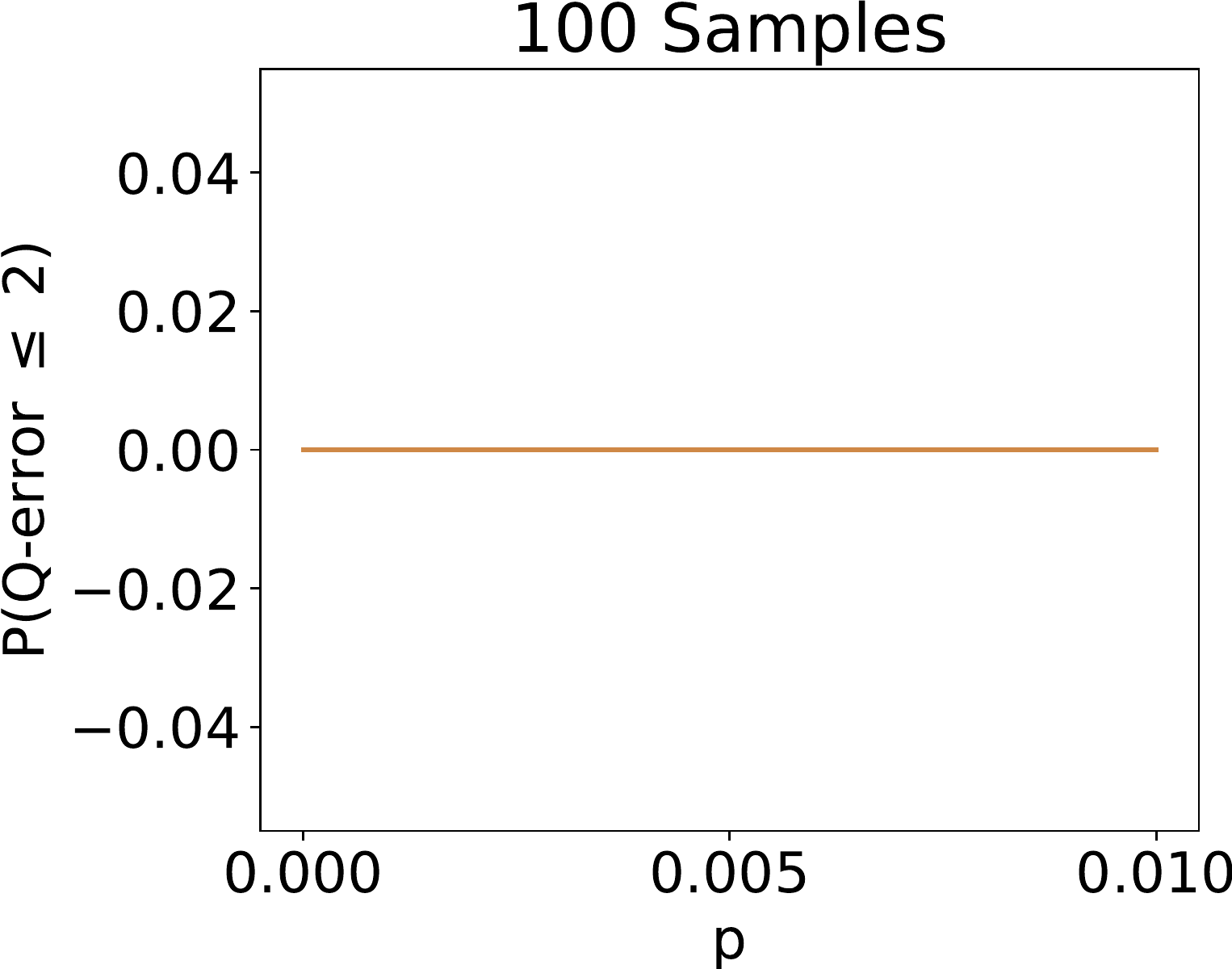}
    \end{subfigure}
    \begin{subfigure}[b]{0.32\textwidth}
        \includegraphics[width=1.01\textwidth]{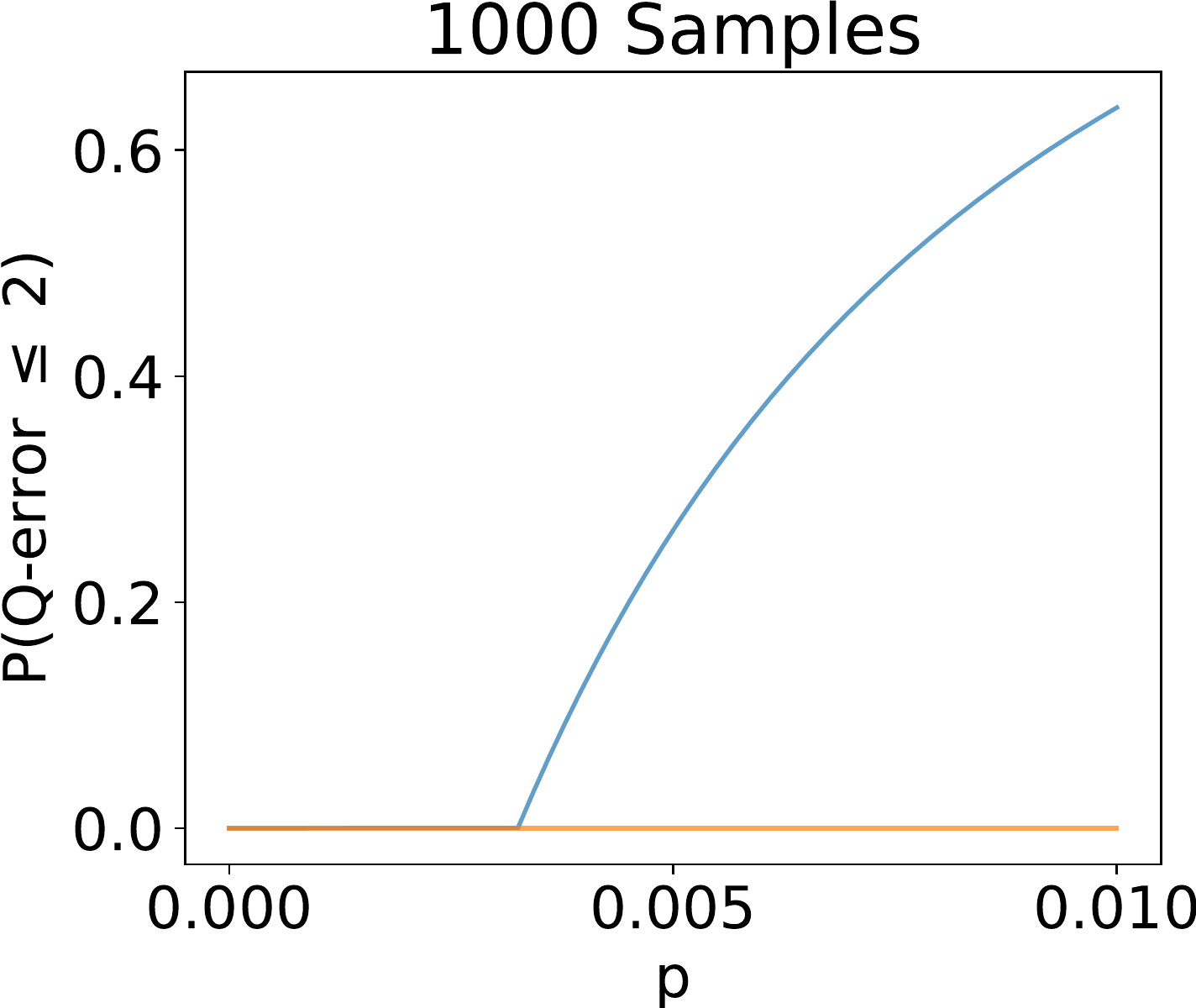}
    \end{subfigure} \vspace{-0.1in}
    \begin{subfigure}[b]{0.60\textwidth}
        \includegraphics[width=1.01\textwidth]{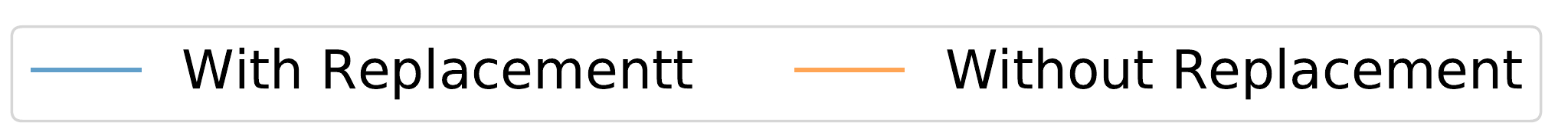}
    \end{subfigure}
    \caption{Comparing random uniform sampling with and without replacement at different cardinality values. Assume the length of the table is 1 billion. }
    \label{fig:compare_with_without}
\end{figure*} 
 
\begin{figure}[htb]
    \centering
     \begin{subfigure}[b]{0.42\textwidth}
        \includegraphics[width=1\textwidth]{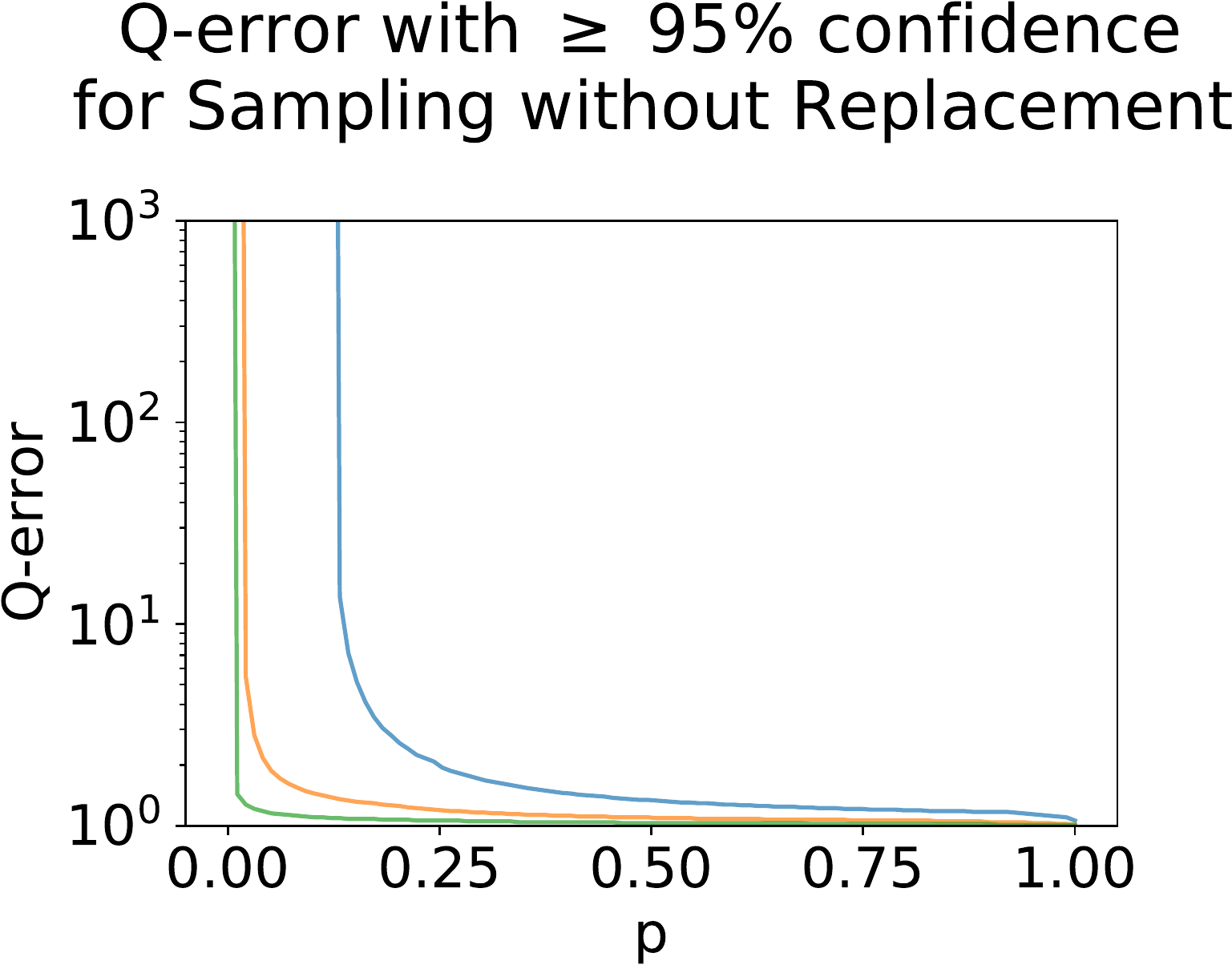}
    \end{subfigure}
     \begin{subfigure}[b]{0.42\textwidth}
        \includegraphics[width=1\textwidth]{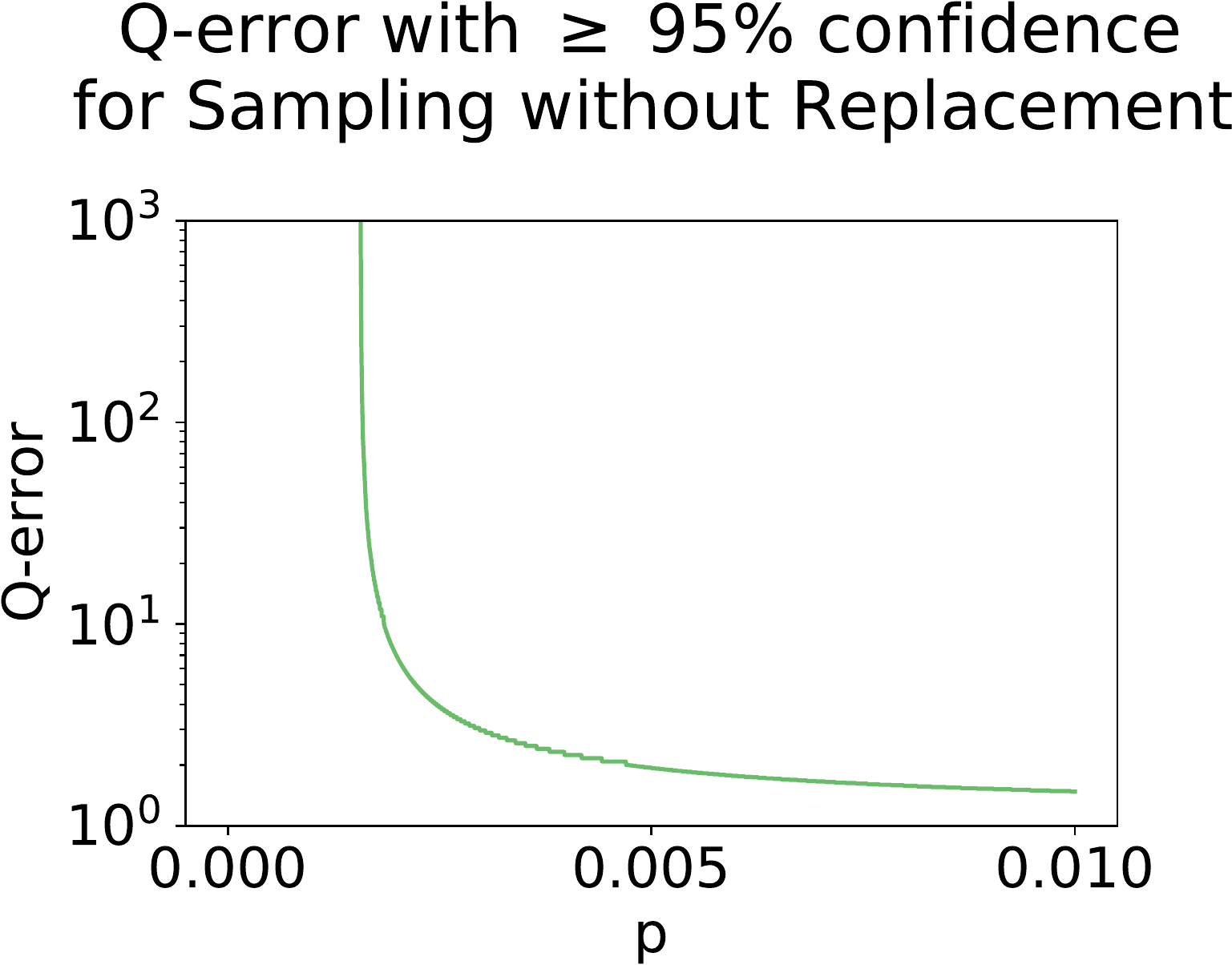}
    \end{subfigure}
    \begin{subfigure}[b]{0.70\textwidth}
        \includegraphics[width=1\textwidth]{fig/legend_n_samples.pdf}
    \end{subfigure} \vspace{-0.1in}
    \caption{Plotting the Q-error with 95\% confidence in terms of $p$.  (Left) Results for $p\in(0,100\%]$.  (Right) Results for $p\in(0,1\%]$.}
    \label{fig:n_samples_norp}
\end{figure} 
 
\begin{figure*}[h!]
    \centering
     \begin{subfigure}[b]{0.32\textwidth}
        \includegraphics[width=1\textwidth]{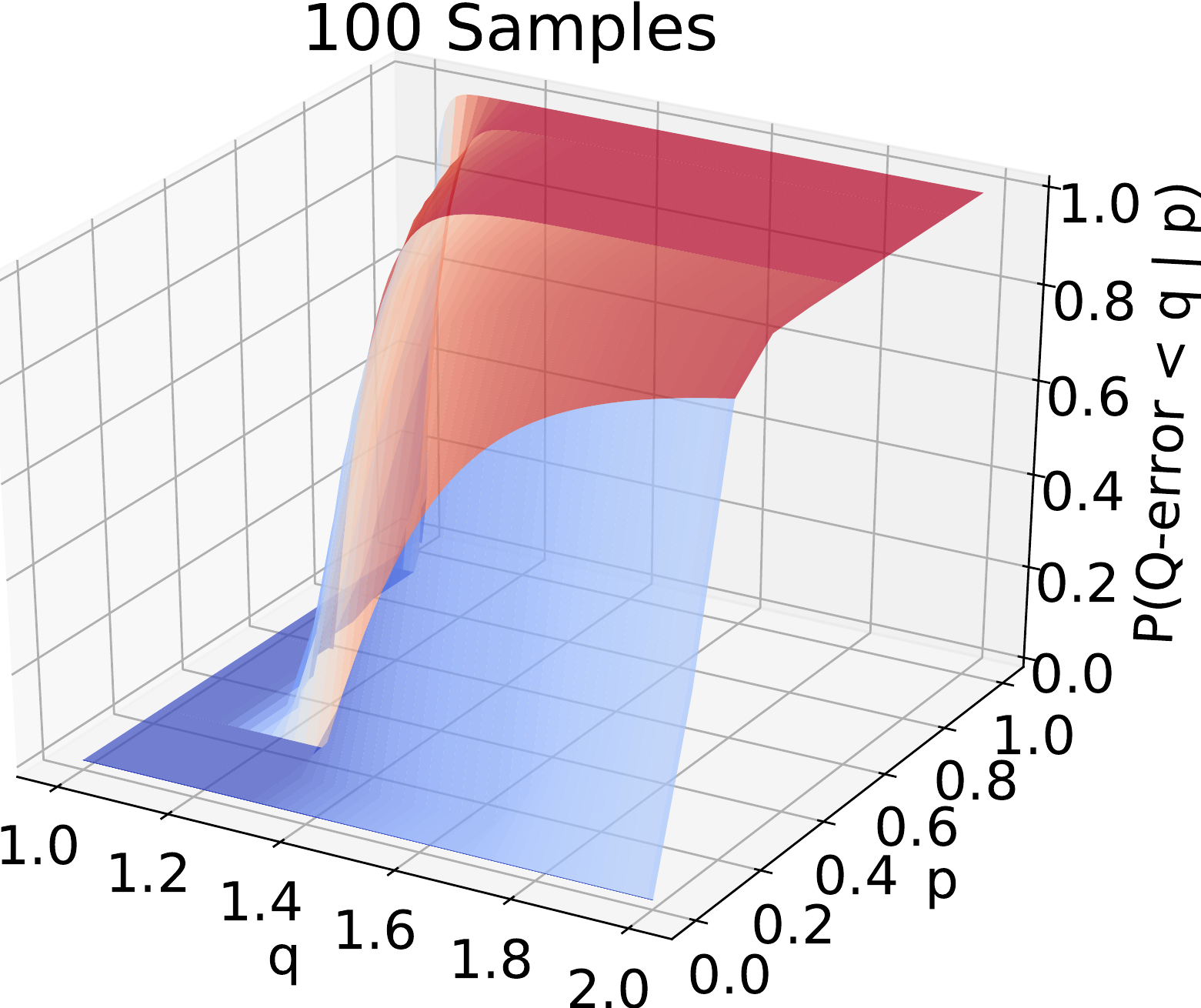}
    \end{subfigure}
    \begin{subfigure}[b]{0.32\textwidth}
        \includegraphics[width=1\textwidth]{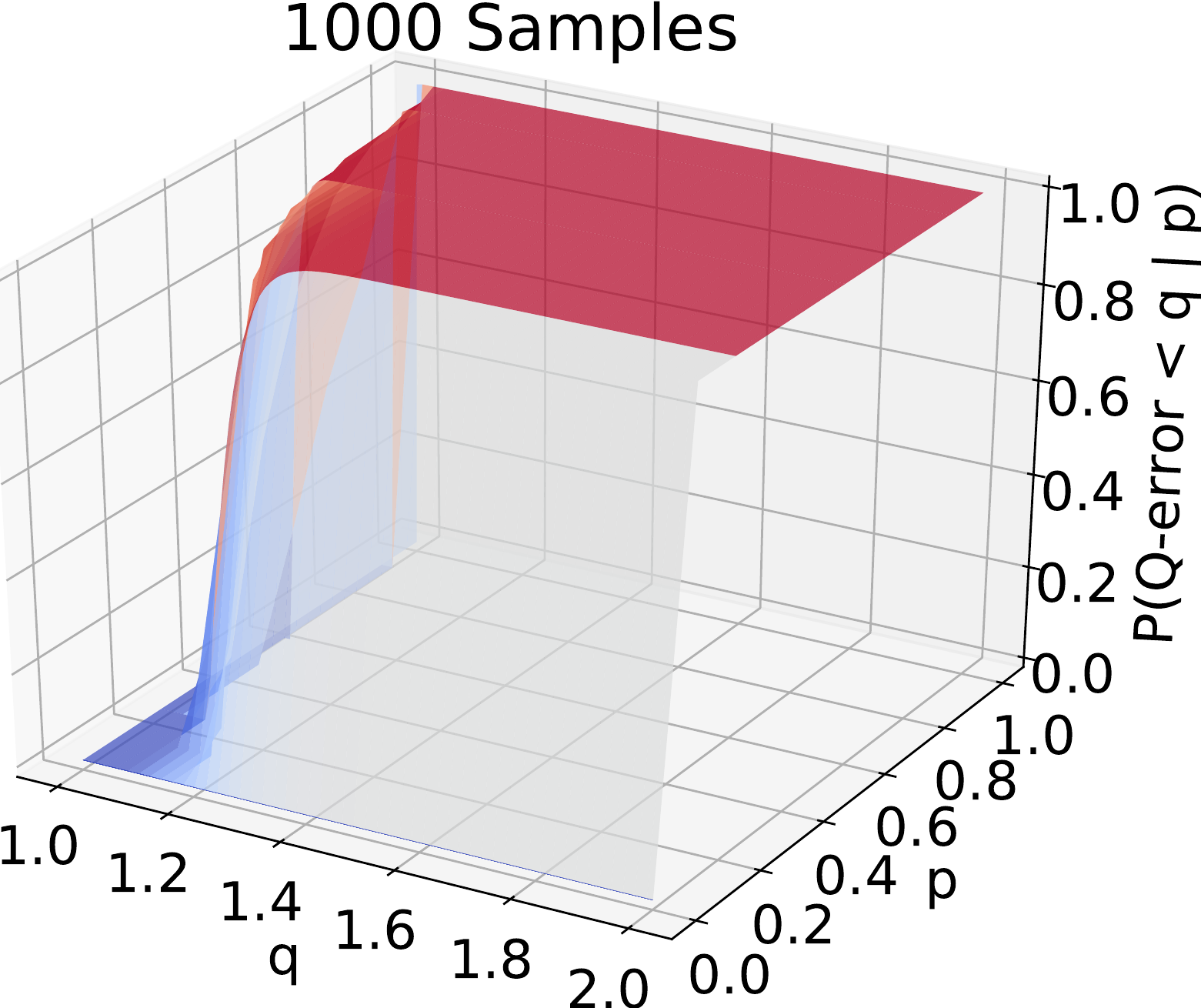}
    \end{subfigure}
    \begin{subfigure}[b]{0.32\textwidth}
        \includegraphics[width=1\textwidth]{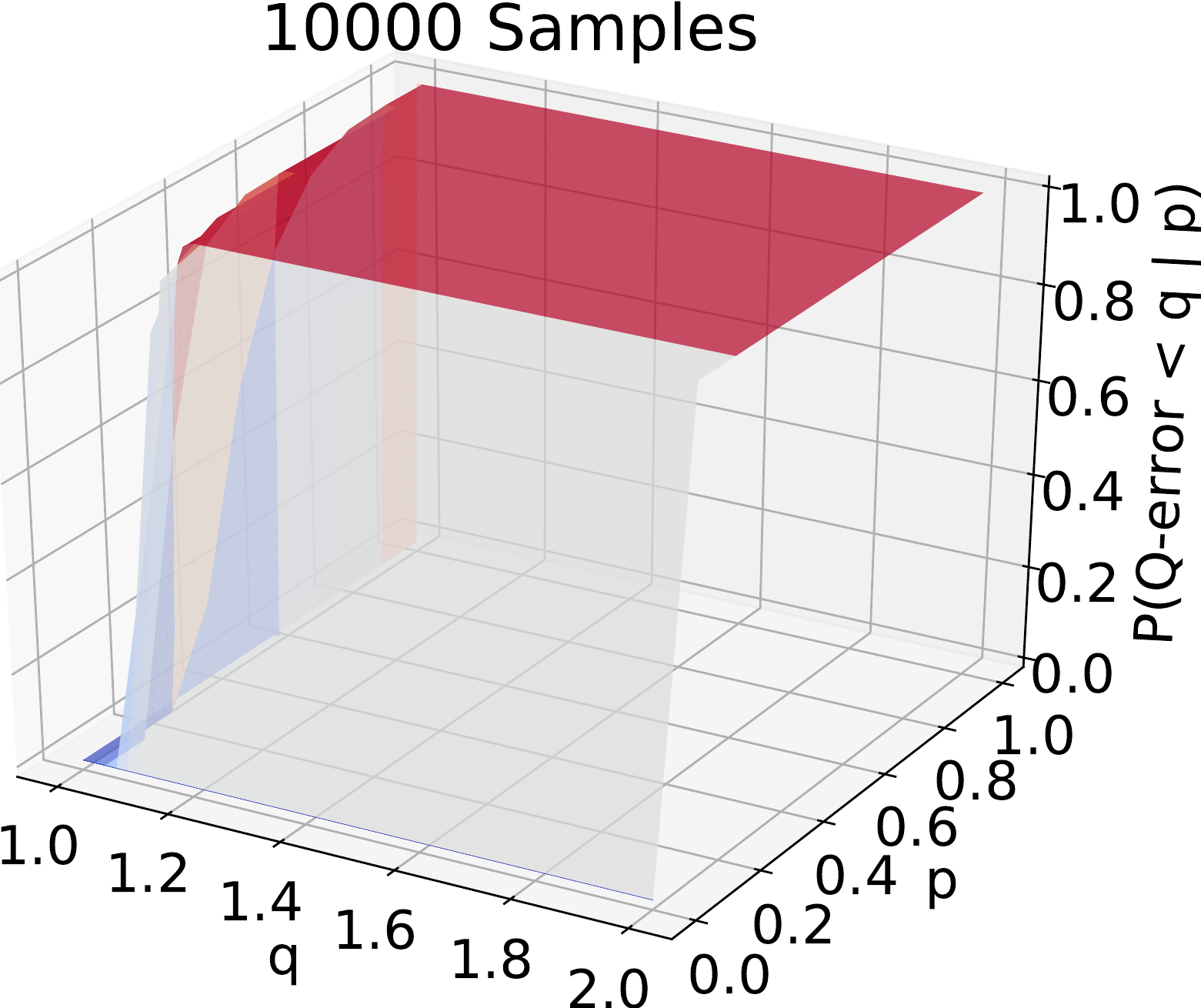}
    \end{subfigure} 
    \caption{3D plotting of the probability that random uniform sampling is better than $q$ at true cardinality $p$. }
    \label{fig:3d_no_rp}
\end{figure*}

The analysis procedure using Bernstein-Serfling inequality is almost identical to the previous steps.
Note there are two roots  for $\delta$, but one root violates the constraint that 
$\cfrac{2 \rho \log(1/\delta)}{k} \geq 0 $.
 
\begin{align*}
\mathds{P}(Pred > q X) 
    &= \mathds{P}\Bigg( \frac{\sum\limits_{t=1}^k (X_t - p)}{k} > pq - p\Bigg) \\
    & \text{(Let $\delta =  \exp \Big(- \frac{k}{\zeta^2} (- \sqrt{ 2 \zeta \rho \sigma^2 (pq - p)  + \rho^2 \sigma^4} + (pq - p)\zeta + \sigma^2 \rho  )  \Big)$)} \\
    &= \mathds{P}\Bigg(\frac{\sum\limits_{t=1}^k (x_t - p)}{k} >  \sigma \sqrt{\cfrac{2 \rho \log(1/\delta)}{k}} + \cfrac{\zeta  \log (1/\delta)}{k}\Bigg)   \\
    &\leq 2 \delta  =  2 \exp \Big(- \frac{k}{\zeta^2} (- \sqrt{ 2 \zeta \rho \sigma^2 (pq - p)  + \rho^2 \sigma^4} + (pq - p)\zeta + \sigma^2 \rho) \Big),
\end{align*}
 
\begin{align*}
    \mathds{P}(Pred < X/q)  
    &= \mathds{P}\Bigg( \frac{\sum\limits_{t=1}^k (Y_t - \bar{p})}{k} > p - p/q \Bigg) \\
    & \text{(  Let $\delta = \exp \Big(- \frac{k}{\zeta^2} (- \sqrt{ 2 \zeta \rho \sigma^2 (p - p/q)  + \rho^2 \sigma^4} + (p - p/q)\zeta + \sigma^2 \rho  )  \Big)$     )}\\
    &= \mathds{P}\Bigg( \frac{\sum\limits_{t=1}^k (Y_t - \bar{p})}{k} >  \sigma \sqrt{\cfrac{2 \rho \log(1/\delta)}{k}} + \cfrac{\zeta \log (1/\delta)}{k}\Bigg)   \\
    &\leq 2 \delta = 2 \exp \Big(- \frac{k}{\zeta^2} (- \sqrt{ 2 \zeta \rho \sigma^2 (p - p/q)  + \rho^2 \sigma^4} + (p - p/q)\zeta + \sigma^2 \rho  )  \Big).
\end{align*}
 
\begin{bound} 
\label{bd:bs}
By applying the Bernstein-Serfling Inequality, the Q-error of random uniform sampling without replacement is bounded by
\begin{align*}
\mathds{P}(\text{Q-error $\leq$ q}) &= 1 - \mathds{P}(\text{Q-error $>$ q}) \\
& \geq 1 
 - 2 \exp \Big(- \frac{k}{\zeta^2} (- \sqrt{ 2 \zeta \rho \sigma^2 (pq - p)  + \rho^2 \sigma^4} + (pq - p)\zeta + \sigma^2 \rho  )  \Big) \\
& - 2 \exp \Big(- \frac{k}{\zeta^2} (- \sqrt{ 2 \zeta \rho \sigma^2 (p - p/q)  + \rho^2 \sigma^4} + (p - p/q)\zeta + \sigma^2 \rho  ) \Big).
\end{align*}
\end{bound}
 
\noindent Putting together the $\Omega$ and $\Psi$ in the above bounds, we derive Theorem \ref{theory:without_replacement}.

\begin{figure*}[h!]
    \centering
     \begin{subfigure}[b]{0.32\textwidth}
        \includegraphics[width=1\textwidth]{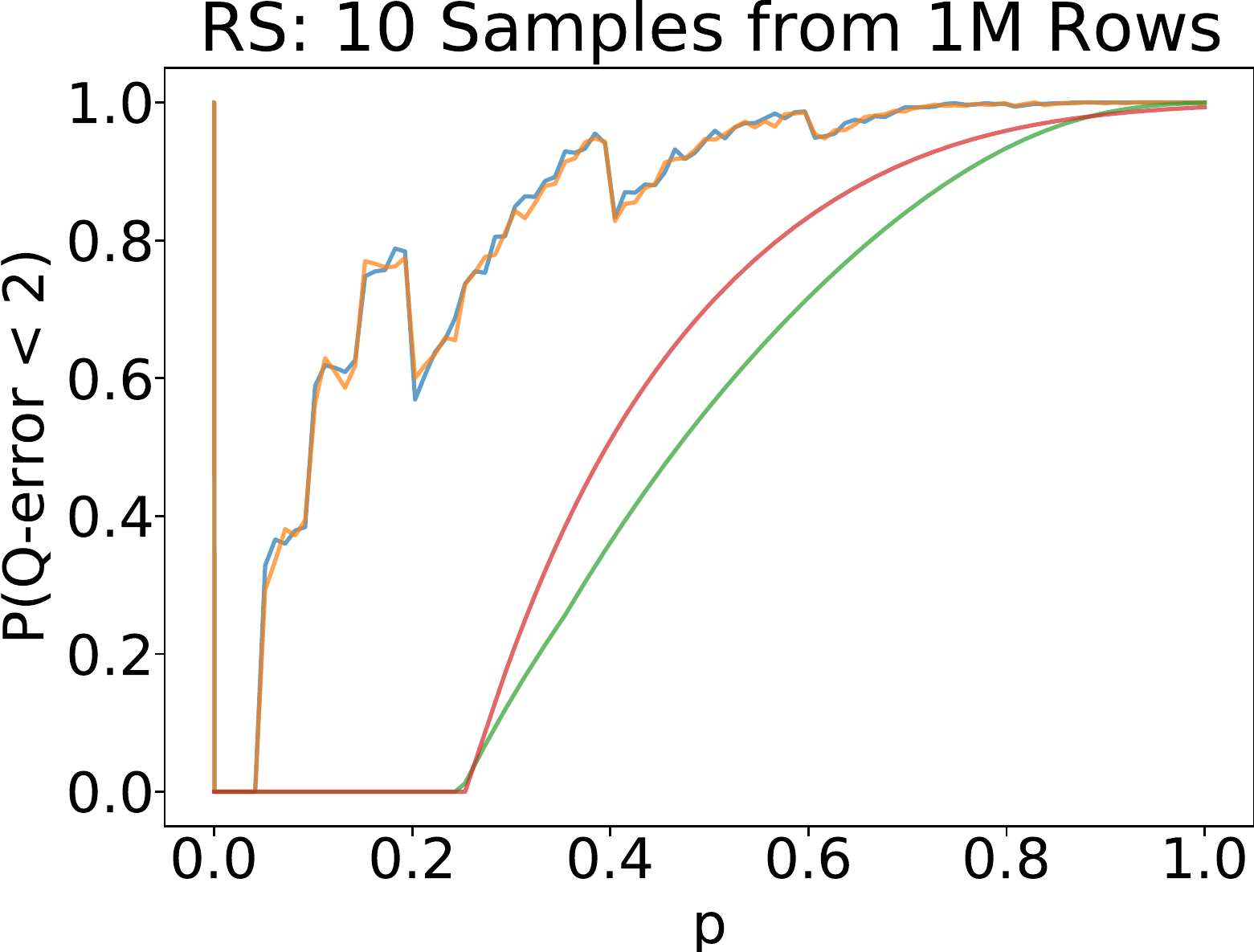}
    \end{subfigure}
    \begin{subfigure}[b]{0.32\textwidth}
        \includegraphics[width=1\textwidth]{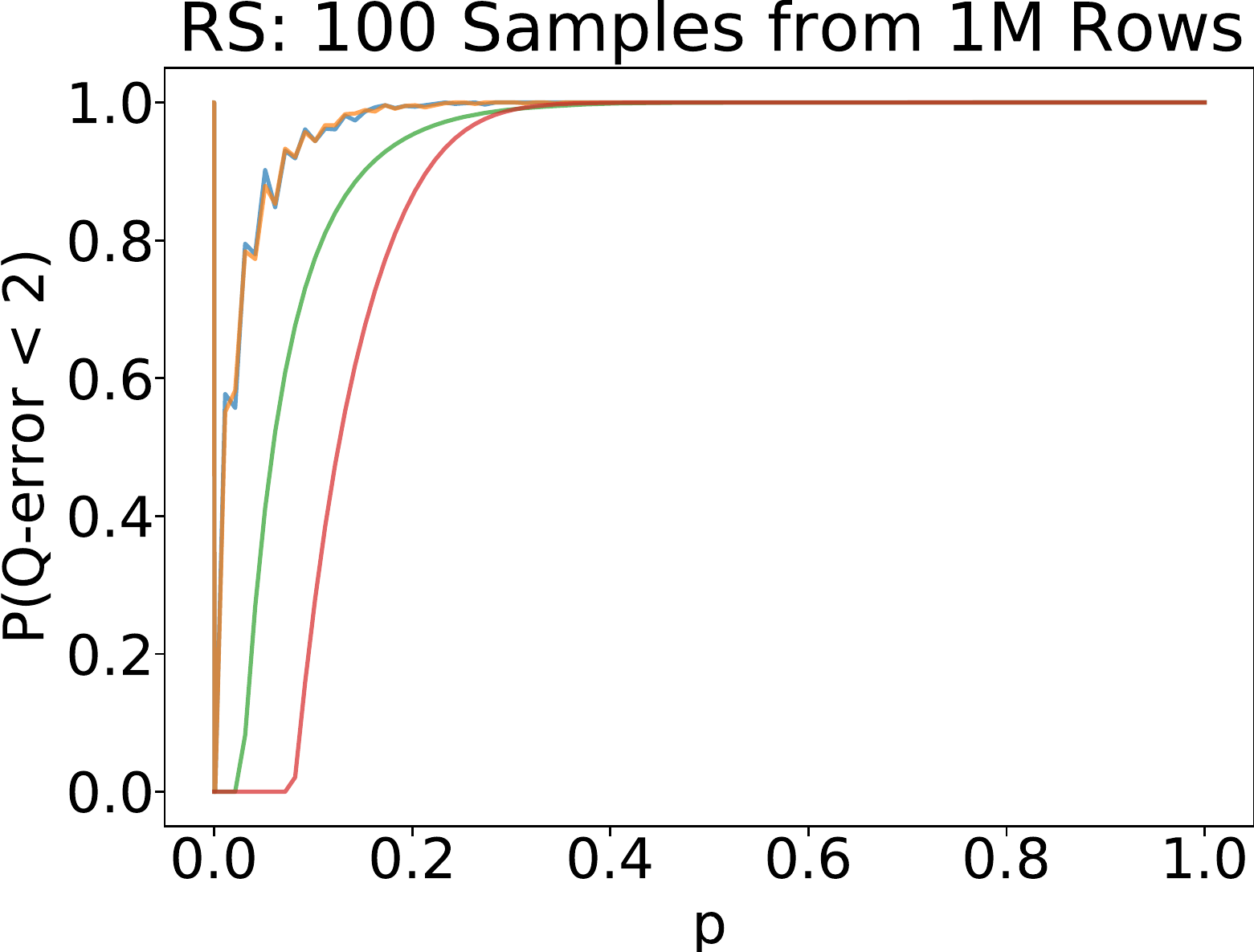}
    \end{subfigure}
    \begin{subfigure}[b]{0.32\textwidth}
        \includegraphics[width=1\textwidth]{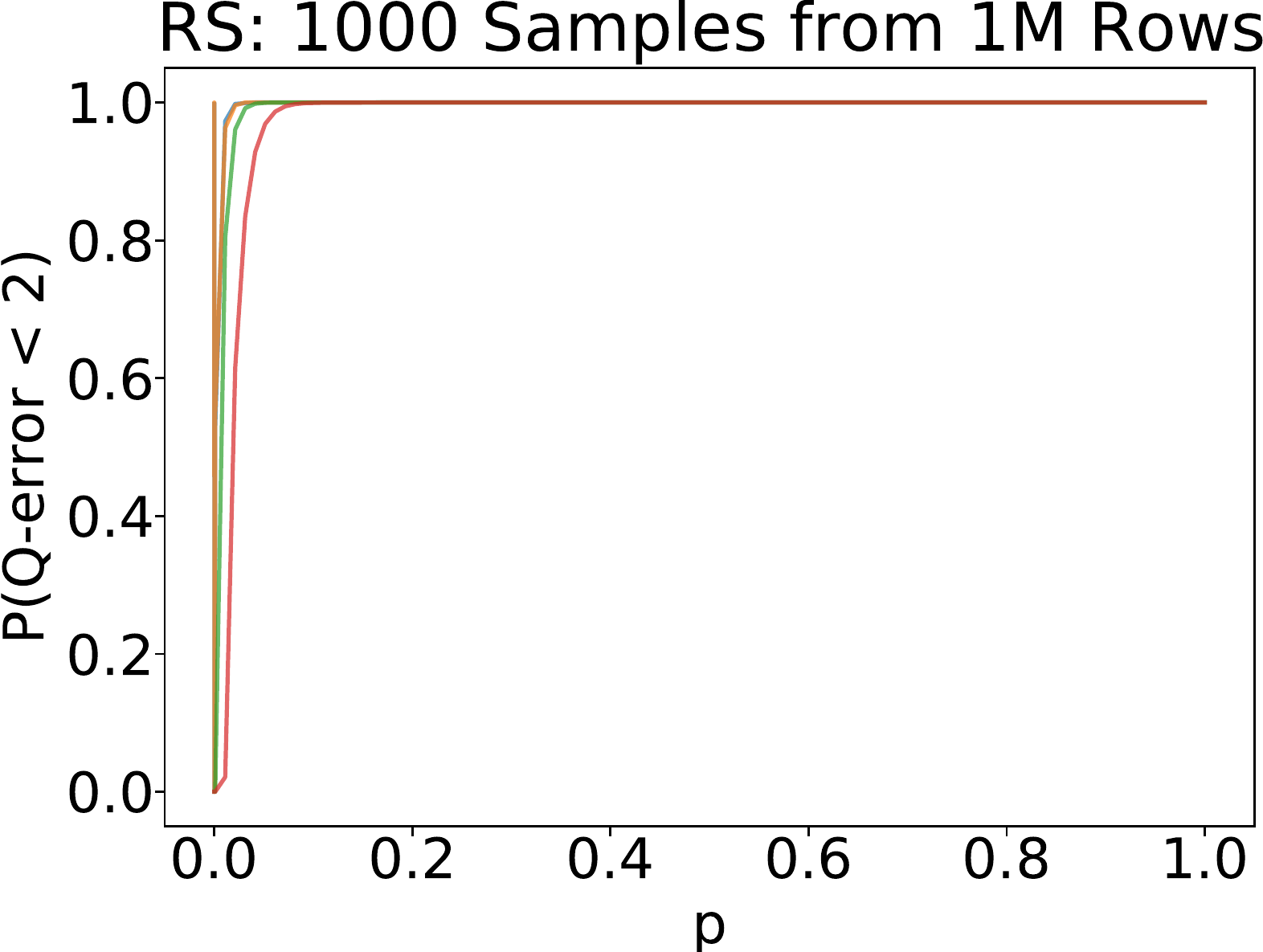}
    \end{subfigure} 
    \begin{subfigure}[b]{0.70\textwidth}
        \includegraphics[width=1\textwidth]{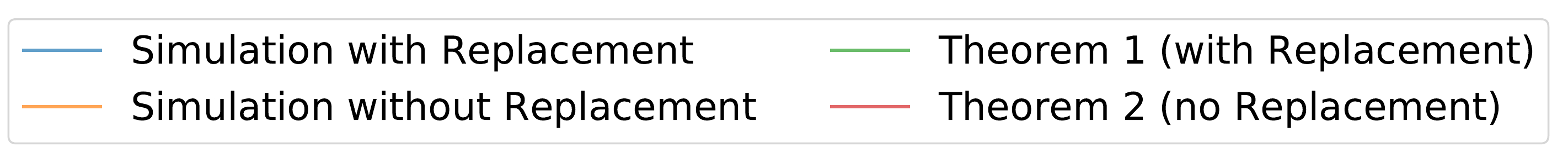}
    \end{subfigure} 
    \caption{For each point in the figure, we conducted a simulation for 1,000 times and plotted the simulation results with Theorem~\ref{theory:with_replacement} and Theorem~\ref{theory:without_replacement} assuming 1 million rows in the table.  When $k$ is small, our bounds are conservative. At a slightly larger $k$, our bounds are tight compared to the simulations.}
    \label{fig:simulation}
\end{figure*}

\vspace{0.1in}\noindent\textbf{Visualization.}
Figure~\ref{fig:hoeffding_bernstein_serfling} plots and compares the bounds in this section. 
With 10K samples, the Q-error is almost always small ($q\le$2) for both bounds. 
Figure~\ref{fig:compare_with_without} compares random uniform sampling with and without replacement at different $p$ values; we observe a tighter bound for random uniform sampling without replacement.  Figure \ref{fig:n_samples_norp} shows the Q-error with at least 95\% confidence at different $p$ values.
Figure~\ref{fig:3d_no_rp} shows the 3-D plotting of the probability that random uniform sampling is better than a given Q-error threshold. 
We also simulate and plot the results in Figure~\ref{fig:simulation}.

\begin{table}[]
\centering
\begin{tabular}{|c|c|c|c|c|c|c|c|}
\hline
\multirow{2}{*}{$p$} & \multirow{2}{*}{C} & \multicolumn{2}{c|}{100 Samples} & \multicolumn{2}{c|}{1000 Samples} & \multicolumn{2}{c|}{10000 Samples} \\ \cline{3-8} 
                       &                    & R       & NR       & R        & NR       & R        & NR        \\ \hline
0.0002 & 166    & 0.00 & 0.00   & 0.00 & 0.00   & 0.00 & 0.00   \\ \hline
0.0003 & 333    & 0.00 & 0.00   & 0.00 & 0.00   & 0.12 & 0.00   \\ \hline
0.0005 & 500    & 0.00 & 0.00   & 0.00 & 0.00   & 0.39 & 0.00   \\ \hline
0.0007 & 666    & 0.00 & 0.00   & 0.00 & 0.00   & 0.56 & 0.00   \\ \hline
0.0008 & 833    & 0.00 & 0.00   & 0.00 & 0.00   & 0.68 & 0.00   \\ \hline
0.0010 & 1000   & 0.00 & 0.00   & 0.00 & 0.00   & 0.76 & 0.00   \\ \hline
0.0017 & 1666   & 0.00 & 0.00   & 0.00 & 0.00   & 0.92 & 0.42   \\ \hline
0.0033 & 3333   & 0.00 & 0.00   & 0.12 & 0.00   & 0.99 & 0.85   \\ \hline
0.0050 & 5000   & 0.00 & 0.00   & 0.39 & 0.00   & 1.00 & 0.96   \\ \hline
0.0067 & 6666   & 0.00 & 0.00   & 0.56 & 0.00   & 1.00 & 0.99   \\ \hline
0.0083 & 8333   & 0.00 & 0.00   & 0.68 & 0.00   & 1.00 & 1.00   \\ \hline
0.0100 & 10000  & 0.00 & 0.00   & 0.76 & 0.00   & 1.00 & 1.00   \\ \hline
0.1667 & 166666     & 0.92 & 0.75   & 1.00 & 1.00   & 1.00 & 1.00   \\ \hline
0.3333 & 333333     & 0.99 & 1.00   & 1.00 & 1.00   & 1.00 & 1.00   \\ \hline
0.5000 & 500000     & 1.00 & 1.00   & 1.00 & 1.00   & 1.00 & 1.00   \\ \hline
0.6667 & 666666     & 1.00 & 1.00   & 1.00 & 1.00   & 1.00 & 1.00   \\ \hline
0.8333 & 833333     & 1.00 & 1.00   & 1.00 & 1.00   & 1.00 & 1.00   \\ \hline
1.0000 & 1000000    & 1.00 & 1.00   & 1.00 & 1.00   & 1.00 & 1.00   \\ \hline
\end{tabular}
  \caption{Confidence that Q-error is at most 2 from Theorem \ref{theory:with_replacement} and Theorem \ref{theory:without_replacement}. Assuming there are 1 million rows, we provide ratio $p$ and the cardinality (C). R = Sampling with Replacement. NR = Sampling without Replacement. $1.00$ is rounded from values that is greater than $0.995$.}
\end{table}
 
\section{Related Work}
\label{sec:related}
 
Recently, \cite{dutt2020efficiently} applied random uniform sampling to dynamically modify the sample size $k$ used in machine learning model training. Their approach can achieve $\epsilon$-optimally by a novel approximation algorithm with Chernoff bound and Hoeffding's inequality. The authors change the sample size $k$ based on the number of satisfied rows in the samples. In this paper, we fix the sample size $k$ and analyze the bound of the Q-error, assuming the ground truth $p$ is known but the number of satisfied rows is unknown. These two analyses may result in different applications.  
In another study, \cite{moerkotte2020alpha} calculated the expected Q-error for TPC-H and several other datasets when the cardinality is extremely small. They also created a novel algorithm that can reduce the expected Q-error. Our analysis in this paper provides confidence intervals other than expected accuracy.

\section{Conclusion}
In this paper, we apply different statistical tools to analyze the upper Q-error bounds of random uniform sampling for single-table cardinality estimation. Our analysis indicates that a simple sampling already provides robust estimates when the true cardinality is relatively high (e.g., at 1000 rows, $> 1\%$ selectivity).

It is easy to see that our analysis in this paper can be extended to sample-based CE after join. Using small samples for each join relation may yield a good error bound; however, deciding which join relations to materialize is an open question. For join on samples, we refer the readers to recent analyses in~\cite{huang2019joins}.

\bibliographystyle{abbrv}
\bibliography{main}

\section*{Appendix}

\begin{theorem}
  \label{theory:with_replacement_with_hoef}
  For cardinality estimation over single tables, 
    the Q-error of random uniform sampling with replacement is bounded by
  \begin{align*}
  & \mathds{P}(\text{Q-error $\leq$ q}) \geq 1 - \Omega - \Psi, where \\
  \Omega &=  \min\Bigg(\Big(\cfrac{e^{q - 1}}{q^{q}} \Big)^{pk}, \exp\Big(- \frac{k (pq-p)^2}{2 \sigma^2 + 2 (pq-p) / 3}\Big), exp\Big(-2p^2(q-1)^2k\Big) \Bigg), \\
    \text{when~} pq > 1,
    \Psi &=  \min\Bigg(\Big(e^{(\frac{1}{q} - 1)} q^{\frac{1}{q}} \Big)^{pk},  \exp\Big(- \frac{k (p-p/q)^2}{2 \sigma^2 + 2 (p-p/q) / 3}\Big), exp\Big(-\frac{2k (pq - 1)^2}{q^2}\Big)  \Bigg), \\
    \text{when~} pq \leq 1,
  \Psi &=  \min\Bigg(\Big(e^{(\frac{1}{q} - 1)} q^{\frac{1}{q}} \Big)^{pk},  \exp\Big(- \frac{k (p-p/q)^2}{2 \sigma^2 + 2 (p-p/q) / 3}\Big) \Bigg).
  \end{align*}
  \end{theorem}
  
The Hoeffding's inequality can be written as the following inequalities, where $\epsilon > 0$:
\begin{align}
\mathds{P}(\hat{X} \geq (p + \epsilon) k) \leq exp(-2 \epsilon^2 k), \\
\mathds{P}(\hat{X} \leq (p - \epsilon) k) \leq exp(-2 \epsilon^2 k).
\end{align}

We can improve Theorem \ref{theory:with_replacement} by adding Hoeffding's inequality.
The bound provided below requires $pq > 1$ for under-estimation and we excluded it from Theorem \ref{theory:with_replacement} for simplicity.

\begin{align*}
\mathds{P}(Pred \geq qX) &= \mathds{P}(\hat{X} \geq q \mu) \\
& \text{(Let $\epsilon = qp - p$)} \\
& \leq exp(-2p^2(q-1)^2k).
\end{align*}

\begin{align*}
\mathds{P}(Pred \leq \frac{1}{q} X) &= \mathds{P}(\hat{X} \leq  \frac{1}{q} \mu) \\
& \text{(Let $\epsilon = p - \frac{1}{q}$ and  $pq > 1$)} \\
& \leq exp(-\frac{2k (pq - 1)^2}{q^2})  \text{~~~~~(if $pq > 1$)}.
\end{align*}

\begin{bound} 
\label{bd:hoeffding}
By Hoeffding's inequality,
when $pq > 1$, 
random sampling with replacement has the Q-error bound
\begin{align*}
\mathds{P}(\text{Q-error $<$ q}) 
 \geq 1 -   exp(-2p^2(q-1)^2k) -  exp(-\frac{2k (pq - 1)^2}{q^2}).
\end{align*}
\end{bound}

\end{document}